\documentclass{amsart}
\usepackage{amsthm,amsmath,amsfonts,amssymb,graphicx,overpic,rotating,url}

\newcommand {\F} {\mathcal{F}}  %
\newcommand {\vst} {v_\text{start}}
\newcommand\EXP{\mathbf{E}}

\newcommand {\defi} {:=}  %
\newcommand {\polar} {\bigtriangleup}
\newcommand {\tr} {\top}

\newcommand {\dto} {\longrightarrow}
\newcommand {\uto} {\text{ \/\textemdash\/ }}
\newcommand {\dtoorfrom} {\longleftrightarrow}

\newcommand {\q} {q}
\newcommand {\s} {s}
\newcommand {\Qset} {Q}
\newcommand {\Sset} {S}
\newcommand {\mqa} {{\bar{\q}}_1}
\newcommand {\mqb} {{\bar{\q}}_2}
\newcommand {\msa} {{\bar{\s}}_1}
\newcommand {\msb} {{\bar{\s}}_2}

\newcommand {\conv} {\operatorname{conv}}
\renewcommand {\exp}[1]{\operatorname{exp}(#1)}

\newcommand{\R}{\mathbb{R}}
\newcommand{\Z}{\mathbb{Z}}
\newcommand{\N}{\mathbb{N}}

\newcommand {\cut} {\operatorname{\delta}}
\newcommand {\incut} {\operatorname{\delta}^{\rm in}}
\newcommand {\outcut} {\operatorname{\delta}^{\rm out}}
\newcommand {\preds} {\operatorname{pred}}
\newcommand {\succs} {\operatorname{succ}}
\newcommand {\length} {\operatorname{length}}

\newcommand {\verts} {\operatorname{vert}}

\theoremstyle{plain}
\newtheorem{lemma}{Lemma}[section]
\newtheorem{theorem}[lemma]{Theorem}
\newtheorem{corollary}[lemma]{Corollary}
\newtheorem{proposition}[lemma]{Proposition}
\theoremstyle{definition}
\newtheorem{definition}[lemma]{Definition}

\begin{document}

  \title{The Random Edge Simplex Algorithm on\\Dual Cyclic 4-Polytopes}
  \author{Rafael Gillmann} %
  \address{MA 6--2\\TU Berlin\\10623~Berlin\\Germany} %
  \email{gillmann@math.tu-berlin.de} %
  \urladdr{http://www.math.tu-berlin.de/~gillmann/} %
  \thanks{%
    The author was supported by the DFG-Forschergruppe \emph{Algorithmen,
      Struktur, Zufall} (FOR 413/1-1, Zi 475/3-1) and by the DFG Leibniz
    grant of G.~M.~Ziegler. %
  } %
  \date{May 4, 2006}%


\begin{abstract}
  The simplex algorithm using the random edge pivot-rule on any
  realization of a dual cyclic 4-polytope with $n$ facets does not
  take more than $O(n)$ pivot-steps.
  This even holds for general abstract objective functions (AOF) /
  acyclic unique sink orientations (AUSO).
  The methods can be used to show analogous results for products of
  two polygons.
  In contrast, we show that the random facet pivot-rule is slow on
  dual cyclic 4-polytopes, i.e. there are AUSOs on which random facet
  takes at least $\Omega(n^2)$ steps.



\end{abstract}

\maketitle

\section{Introduction}
\label{sect:intro}

Linear Programming (LP) is the problem of minimizing some linear
function $x \mapsto c^t x$ in $d$ variables subject to $n$ linear
inequalities. In geometric terms we are given a polyhedron $P$ in
$\R^d$ defined as the intersection of $n$ half-spaces; the objective
is to find some extremal point in $P$ w.r.t.\ a given linear function
$c^tx$.

The \emph{simplex algorithm} is the oldest linear programming
algorithm.  It was devised by Dantzig in 1947 and first published in
1951 in \cite{Dantzig51}.  In terms of geometry it finds the minimal
vertex of the given simple $d$-Polytope $P$ by starting at a given
starting vertex and iteratively moving to an improving neighbor until
the minimal vertex is reached.  Usually there are several improving
neighbors to choose from. A \emph{pivoting rule} decides which one to
pick. So the simplex algorithm is actually a class of algorithms and
we will refer to the simplex algorithm with a certain pivot-rule also
by the name of the pivot-rule.

In 1980 Khachiyan proved that LP can be solved in polynomial time by
the ellipsoid method \cite{Khachiyan80} depending on $d$, $n$ and the
bit-size of the input (Turing-Machine model).
Until now there is no (combinatorial) strongly polynomial algorithm
known to solve LP in running time bounded by a polynomial depending on
$d$ and $n$ only (unit-cost model / RAM model). The simplex algorithm
seems to be a natural candidate for such an algorithm, since almost
every reasonable pivot-rule chooses the next vertex in strongly
polynomial time. Thus the running time can be expressed as the number
of pivoting steps.
There is no pivot-rule known which requires only a polynomial number
of pivoting steps.

For many pivot-rules difficult inputs have been constructed on which
an exponential number of pivot-steps are required (in the dimension
$d$).  The first examples were the famous Klee-Minty Cubes due to Klee
and Minty in 1972 \cite{KM72} which showed that Dantzig's original
pivot-rule \cite{Dantzig63} could visit all vertices in a cube and
thus requires an exponential number of steps. In fact for most
deterministic pivot-rules such examples are known, c.f. the overview
by Goldfarb 1994 \cite{Goldfarb94}. Many of these constructions have
been unified by Amenta and Ziegler's deformed products \cite{AZ96}.

Two strategies have mainly been followed to try to overcome the
exponential worst-case behavior of the simplex algorithm.
The first idea is to investigate the average case rather then
the worst case.  Borgwardt showed in 1987 that over random LPs
w.r.t. a certain probability distribution the shadow-vertex
simplex-algorithm needs only polynomial many pivot-steps
\cite{Borgwardt87}.
In \cite{ST04} Spielman and Teng introduced the smoothed analysis
which combines advantages of worst-case and average-case analysis.
Smoothed analysis measures the maximum of the expected running time
over inputs under small random perturbations.  They prove that the
simplex algorithm with the shadow-vertex pivot rule has a polynomial
smoothed complexity, i.e. the running time is polynomial in the input
size and the standard deviation of Gaussian perturbations.
Recently Kelner and Spielman \cite{KS05} introduced a randomized
``simplex like'' algorithm which runs in polynomial (but not strongly
polynomial) time.  Their algorithm solves a randomized sequence of LPs
using the shadow-vertex simplex algorithm building on \cite{ST04}.

The second idea is to introduce randomness to the pivoting rule and
not the input yielding \emph{randomized pivot-rules}.  In the
following we will consider the worst-case expected running-time of
randomized pivot-rules.
The first substantial progress on upperbounds on randomized pivot
rules was obtained by Kalai in \cite{Kalai92} and independently by
Matou\v{s}ek, Sharir, and Welzl in \cite{MSW96}. Kalai proved that
random facet needs at most (in expectation) $\exp{ O(\sqrt{d \log d})
}$ steps. This was the first sub-exponential running time for any
pivot rule. Matou\v{s}ek, Sharir, and Welzl had a similar result.

The analysis of random facet relies on rather simple and general
properties of orientations of polytope graphs induced by linear
objective functions.  These properties establish a more general purely
combinatorial framework in which the proof works.  In this paper we
will call those orientations %
\emph{acyclic unique sink orientations (AUSO)} %
though Kalai introduced them as abstract objective functions (AOF) in
\cite{Kalai88}.  AOFs and AUSOs are essentially the same.
The same concept was also introduced independently by Williamson Hoke
in \cite{WH88}.
Besides AUSOs there are also other abstract settings like Sharir and
Welzl's LP-type problems \cite{SW92} and G\"artner's abstract
optimization problems \cite{Gartner95}.

The upper bounds on random facet established in \cite{Kalai92} and
\cite{MSW96} are nearly tight in the setting of AUSOs. In
\cite{Matousek94} Matou\v{s}ek constructs a family of abstract cubes
(AUSOs on cubes) such that random facet requires
$\exp{\Omega(\sqrt{d})}$ steps.  So geometry must help to get under
the sub-exponential bound.  G\"artner showed in \cite{Gartner98} that
on the realizable examples of \cite{Matousek94} 
random facet needs $O(d^2)$ steps only.  An AUSO is called realizable
if there exist an embedded polytope and a linear function such that
the orientation induced by the linear function on the polytope's graph
is the given AUSO.

\paragraph{Random Edge}
The \emph{random edge rule} is probably the most straight forward
randomized pivot-rule: ``Choose the next vertex uniformly at random
among all improving neighbors.''  The price for its simplicity is that
it does not use any polytope specific combinatorial or geometric
information.  Thus it seems reasonable that obtaining good upper
bounds on random edge might be more difficult.
In fact it is already quite difficult to analyze random edge on
3-polytopes.  On 3-polytopes all pivot rules need at most linearly
many steps.  In \cite{KMSZ05} Kaibel, Mechtel, Sharir, and Ziegler
compute the coefficients of linearity of various pivot rules and
random edge turns out to be the most difficult to analyze.

Broder et al.\ showed in \cite{Betal95} that random edge can be
exponential in the height. The height is the shortest directed path
from the unique minimal to the unique maximal vertex of the given
polytope w.r.t.\ the given linear objective function.
On Klee-Minty cubes random edge needs $\Theta(d^2)$ steps only. This
is a result of Balogh and Pemantle \cite{BP05} improving an
earlier result of G\"artner, Henk, and Ziegler \cite{GHZ98}.
G\"artner et al.\ in \cite{Getal03} analyzed random edge on
$d$-polytopes with $d+2$ facets--that is one facet more than the
simplex.
It can be shown that on the abstract cubes from \cite{Matousek94}
random edge only needs $O(d^2)$ steps.
In a survey article by Kalai from 2001 \cite{Kalai01} random edge is
the first among six pivot rules suggested for deeper study.

Up to quite recently the hope was that random edge could be quadratic,
e.g. in $O(dn)$.  But this hope was partially destroyed by
Matou\v{s}ek and Szab\'{o} \cite{MS04} who constructed a family of
abstract cubes on which random edge would need at least %
$\exp{\Omega(d^{1/3})}$ %
steps\footnote{They believe their analysis could be sharpened to
  $\exp{\Omega(d^{1/2})}$.}  %
with high probability.  Thus random edge is exponential.\footnote{It
  would be better to use the term ``super-polynomial'' since the same
  bound--as an upper bound--is considered to be sub-exponential.} %
It seems reasonable to believe that these abstract cubes are not
realizable with high probability.

A trivial general upper bound is the maximal number of vertices of any
$d$-polytope with $n$ facets.  This number is given by the Upper Bound
Theorem (cf. \cite{Ziegler98}) as the number of vertices of the dual
cyclic $d$-polytopes.  G\"artner and Kaibel gave the first non-trivial
general upper bound of $O(N/\sqrt{d})$ in \cite{GK05-pre}, where $N$
denotes the number of vertices of the given polytope.
Thus in contrast to most exponential examples for deterministic pivot
rules, random edge skips a substantial amount of vertices.

A substantial progress in the study of random edge would be a
sub-exponential upper bound similar to the one for random facet.  It
would be great to have a polynomial upper bound but then geometry must
help due to Matou\v{s}ek and Szab\'{o}'s abstract cubes.  

\paragraph{4-polytopes}
The study of random edge on small problems can reveal interesting
properties to attack the general problem or to give insights which
lead to new applicable geometric properties.
In dimension $d=3$ every simple $3$-polytope has exactly $2n-4$ many
vertices.  Thus every pivot rule is linear.  Also it is easy to
construct arbitrary 3-polytopes due to Steinitz' Theorem which
characterizes all graphs of 3-polytopes (see \cite{Ziegler98} for
details).
                                
Dimension $d=4$ is more interesting as $4$-polytopes can have
quadratically many vertices.  $4$-polytopes also admit a richer and
more difficult structure than 3-polytopes.  There is no such theorem
as Steinitz' Theorem known for 4-polytopes.  In fact we do not know
many constructions of simple 4-polytopes with many vertices such that
we explicitly know their combinatorics.

As random edge is a completely combinatorial pivot rule we do not need
any geometric information about the polytope such as coordinates,
angles, or objective function values.  
Two polytopes $P$ and $Q$ are called \emph{combinatorial equivalent}
if they have the same combinatorial structure, i.e. if there exists a
bijection between their face lattices (see~\cite{Ziegler98} for
details).
                                
A very simple example of a 4-polytope is the product of two polygons
with $k$ and $\ell$ vertices.  These polytopes have $k+\ell$ facets
and $k \ell$ vertices.  The combinatorics are very easy and already
many techniques developed in this paper apply for this example.
A more complicated example are the dual cyclic 4-polytopes which have
the most vertices for a given number of facets.  The combinatorics are
known due to Gale's Evenness Condition (cf.  Theorem~\ref{thm:gale}).

\subsection {Results}

The main result of this paper is that random edge is fast on dual
cyclic 4-poly\-topes.  Besides the main result, there are many other
results and ideas presented as part of the proof.
Let us first state the main result in more details.  
\begin{theorem} \label{thm:main}
  Let $P$ be a 4-polytope with $n$ facets which is combinatorially
  equivalent to the dual of the cyclic 4-polytope with $n$ vertices,
  which is defined by taking the convex hull of $n$ points on the
  moment curve. Given any AUSO of $P$ and an arbitrary vertex $\vst$
  on $P$, the random edge simplex algorithm starting at $\vst$ on the
  AUSO of $P$ takes at most $O(n)$ pivot steps in expectation to reach
  the minimal vertex (global sink) of $P$.
\end{theorem}
Since every linear orientation of a polytope $P$ is an AUSO of $P$, a
special case of the Main Theorem is obtained by replacing AUSO with
linear orientations defined by a given linear function on $P$.

Section \ref{sect:analysis} is devoted to the proof of the main
result.  Section \ref{sect:prelims} gives the ideas along which the
proof proceeds.  It analyzes the combinatorial structure and its
interplay with AUSOs on dual cyclic 4-poly\-topes, and provides a
technical framework which is used in the proof of the main result.
Though Section \ref{sect:prelims} analyzes the duals of cyclic
4-polytopes only, the results can be obtained for other 4-polytopes as
well, e.g. on products of two polygons, which is sketched in
Section~\ref{sect:polygons}.

In Section~\ref{sect:rf} we sketch the construction of AUSOs on dual
cyclic 4-poly\-topes such that random facet started at the global source
takes at least $\Omega(n^2)$ pivot-steps.
Thus random edge is faster on dual cyclic 4-polytopes than random
facet.

The following is an outline and a short summary of the upcoming
sections devoted to the proof of Theorem~\ref{thm:main}.
Throughout $n$ denotes the number of facets of the dual cyclic
4-polytope considered.

\paragraph{How to Bound the Running Time}
It is easy to deduce an explicit recursion formula to compute the
running time for each vertex for random edge.  
Another approach is to look at the random path defined by starting
random edge at a specific starting vertex.  Such a random path is
equivalent to a flow sending one unit of flow from the starting vertex
to the global sink (minimal vertex) such that at each vertex the
inflow is uniformly distributed over all out-edges. The exact running
time for a specific starting vertex can thus be computed as the costs
of this flow, i.e. the sum over all edge probabilities.
Both was done in \cite{KMSZ05} to derive lower and upper bounds on the
running time's coefficient of linearity for 3-polytopes.  And it turns
out to be quite tedious (already for 3-polytopes).

Let $G=(V,A)$ be a directed graph (with maximal degree at most four)
and $f:V \to \Z$ a \emph{monotone decreasing function}, i.e. $f$ is
not increasing along any directed arc in G.  If every vertex has a
decreasing direct successor w.r.t. $f$, we call $f$ \emph{effectively
  decreasing}.  Then we can bound the running time of random edge by
$O(\# f(V))$.  The monotonicity of $f$ guarantees that we will never
revisit a set of vertices with equal $f$-value after increasing the
$f$-value.  Effectively decreasing ensures that the $f$-value is
increased with probability at least $\tfrac 14$ in every step of
random edge.

Of course it can be difficult to find such a function $f$ for the
whole graph.  But it is already enough to find such functions for each
vertex set in a vertex partition $\Pi$, as long as the quotient graph
$G / \Pi$ (obtained by identifying all vertices in each $W \in \Pi$
and removing double arcs and loops) is still acyclic and the size of
$\Pi$ is bounded by a constant.
In order to find a suitable partition $\Pi$ and functions $f$, we will
need a few more concepts.

\paragraph {Combinatorics of Dual Cyclic 4-Polytopes}
The combinatorics of cyclic polytopes and thus of their duals is
completely known due to Gale's Evenness Condition (c.f. Theorem
\ref{thm:gale}).  There are two kinds of 2-faces: small ones
(triangles and quadrangles) and large ones ($(n-2)$-gons).  The large
2-faces are most interesting for us and are denoted by $\F$.  Each is
a separating cycle of the graph.  They altogether cover all vertices
and edges of the graph.  Moreover they come with a natural
neighborhood relation where two large 2-faces $F, G \in \F$ are
neighbors if their vertices $\verts(F) \cup \verts(G)$ are the
vertices $\verts(H)$ of a 3-face (facet) $H$.  The 2-faces in $\F$ are
numbered $F_0, F_1, \ldots, F_{n-1}$ such that two 2-faces $F_i$ and
$F_j$ are neighbors if and only if $(i-j) \equiv \pm 1 \mod n$.

Gale's Evenness Condition also allows us to draw nice pictures of the
graphs of dual cyclic 4-polytopes. Figure \ref{fig:dcp-graph} depicts
such a graph in a way which illustrates the combinatorics quite well.

\paragraph {AUSOs on Dual Cyclic 4-Polytopes}
Each 2-face $F_i$ has a unique sink and source which we will denote by
$\s_i$ and $\q_i$ respectively. Consider two neighbors $F_i$ and
$F_{i+1}$. Since their vertices span a 3-face, either $\q_i$ or
$\q_{i+1}$ is the source of the 3-face. And thus there is a path
$\gamma_{i,i+1} : \q_i \dtoorfrom \q_{i+1}$ either directed from $\q_i$
to $\q_{i+1}$ or vice versa.  Iterating this construction results in an
(undirected) cycle passing through all sources of the 2-faces in $\F$.
The cycle has at most two sinks, denoted by $\mqa$ and $\mqb$. 
And of course we can apply the whole construction to the sinks of the
2-faces in $\F$ yielding a cycle of sinks which has at most two
sources denoted by $\msa$ and $\msb$.

\paragraph {Intersecting Paths}
The last ingredient for the proof is a simple consequence of the
Jordan Curve Theorem.  The interesting parts of the graphs we consider
a actually planar.  That means that certain paths must intersect in a
certain way.  We introduce the abstract framework of fences.
They allow us to apply the same results to a wider range of polytopes,
e.g. for products of two polygons.

\paragraph {Proving Theorem \ref{thm:main}}
Finally we will put all these ingredients together to proof the main
result.  The idea is to split the set of vertices $V = \verts(P)$ of
$P$ into a constant number of vertex sets. For each vertex set we
define a function.  Now it remains to show that each of the functions
is effectively decreasing. To prove the latter, we will use the cycle
of sources respectively sinks defined earlier and the conditions of
intersecting paths.

\subsection{Notation}

For an arbitrary polytope $P$ we denote by $G(P)$ the (undirected)
graph of the polytope, i.e. its 1-skeleton.  The polytopes considered 
are 4-dimen\-sional and simple, i.e. each vertex is incident to exactly
four facets.  Thus in $G(P)$ each vertex has degree exactly four.

In general an undirected graph $G=(V,E)$ is given by the vertex set
$V$ and the set of edges $E \subset \{ \{x,y\} : x,y \in V, x \not= y
\}$. A digraph $D=(V,A)$ is defined by the vertex set $V$ and the set
of arcs $A \subset V \times V$. All graphs that we consider are
simple, i.e.  there are no parallel edges, arcs or loops.
                                
For any two subsets $V_0, V_1 \subset V$, $E(V_0,V_1)$ denotes
the edges between $V_0$ and $V_1$:
\begin{equation*}
  E(V_0,V_1) = \begin{cases}
    \left\{ \{v,w\} \in E : v \in V_0, w \in V_1 \right\}&
    \text{undirected graph} \\
    \left\{ (v,w) \in A : v \in V_0, w \in V_1 \text{ or } w \in V_0
      , v \in V_1 \right\} &\text{directed graph}
  \end{cases}
\end{equation*}
We define the following abbreviated notations (cut, out-cut, and
in-cut): 
\begin{align*}
  \cut(W) &\defi E(W,V \setminus W) \cup E(V \setminus W, W)\\
  \outcut(W) &\defi \{ (v,w) \in A \;:\; v \in W \text{ and } w \in
  V\setminus W \}\\
  \incut(W) &\defi \{ (v,w) \in A \;:\; v \in V \setminus W \text{ and
  } w \in W \} 
\end{align*}

In the following we work a lot with paths and cycles.
A path of length $\ell$ is a sequence of vertices $v_0, v_1, \ldots,
v_\ell$, s.t. $\{v_i, v_{i+1}\} \in E$ or in case of digraphs
$(v_i,v_{i+1}) \in A$ or $(v_{i+1},v_i) \in A$.  This implies that
$v_i \not= v_{i+1}$.  A path is directed if either all $(v_i,v_{i+1})
\in A$ or all $(v_{i+1},v_i) \in A$.  By $-[v_0,v_1, \ldots, v_\ell]$
we denote the path $[v_\ell, \ldots, v_1, v_0]$.  We use the following
notation for paths:
\begin{itemize}
\item $\gamma : v \uto w$ undirected path connecting $v$ and $w$.
\item $\gamma : v \dto w$ directed path from $v$ to $w$.
\item $\gamma : v \dtoorfrom w$ directed path from $v$ to $w$ or from
  $w$ to $v$.
\end{itemize}

The empty path is a path with just a single vertex and no edge.  Two
paths can be concatenated: 
\begin{equation*}
  [v_0, v_1, \ldots, v_\ell] \circ [w_0, w_1,
\ldots, w_m] \defi [v_1,\ldots, v_\ell, w_2, \ldots, w_m], 
\end{equation*}
$v_\ell = w_0$ must hold!
Let $\gamma \defi [v_0, v_1, \ldots, v_\ell]$ be a path. Then
\begin{equation*}
  \gamma|_{[v_i,v_k]} \defi [v_i, v_{i+1}, \ldots, v_{k-1}, v_k]
\end{equation*}
with $0 \leq i < k \leq \ell$ denotes a sub path of $\gamma$.

A cycle (undirected or directed) is an (undirected resp. directed)
path with $v_0 = v_\ell$.

Note that for a vertex subset $W \subset V$ and a path $\gamma$ the
intersection $\gamma \cap W$ is defined as the set of all vertices of
$\gamma$ in $W$. 

With the notion of directed paths, we can easily define the
predecessors and successors of a given vertex $v$ in a digraph.
\begin{align*}
  \preds(v) \defi& \{ w \in V : \exists \gamma : w \dto v \text{,
    $\gamma$ not empty} \}\\
  \preds[v] \defi& \{ w \in V : \exists \gamma : w \dto v \}\\
  \succs(v) \defi& \{ w \in V : \exists \gamma : v \dto w \text{,
    $\gamma$ not empty} \}\\
  \succs[v] \defi& \{ w \in V : \exists \gamma : v \dto w \}
\end{align*}
Thus $\preds[v] = \preds(v) \cup \{v\}$ and $\succs[v] = \succs(v)
\cup \{v\}$.

A random path from $v \in V$ to $w \in V$ is a function $p : A \to
\R_+$ such that $p$ is a $v$-$w$-flow of value $1$ in the network
$(D,u,v,w)$.  The capacities $u(e)$ are not necessary and are set to
$u(e) = \infty$ for all edges $e \in A$.  Thus $p$ assigns a
probability to each arc.  The random paths defined by random edge are
exactly those from $v$ to the global sink $s$, where for all vertices
$w \in V$ all out-edges $e, e' \in \outcut(w)$ have the same
probability $p(e) = p(e') = {1}/{\lvert \outcut(w)\rvert}$.  The
expected length of $p$ is $\EXP[p] = \sum_{e \in A} p(e)$.

\section {Preliminaries}
\label{sect:prelims}

\subsection{Bounding the Running Time Using Monotone Functions}
\label{sect:prelims:runtime}

This section covers the main ideas of the upcoming runtime analysis of
random edge.
\begin{definition} \label{def:dds} Given a directed graph $D=(V,A)$,
  $W \subset V$, a function $\lambda : V \to \Z$ is called
  \emph{monotone decreasing} on $D$ if $\lambda(v) \geq \lambda(w)$
  for all $(v,w) \in A$. 
  The function $\lambda$ is \emph{effectively decreasing(with respect
    to $W$)} if it is monotone decreasing and for every $v \in V$ that
  is not a global sink, there is a $(v,w) \in A$ such that $\lambda(v)
  > \lambda(w)$ or $w \notin W$, i.e. for every $v$ there is a
  \emph{decreasing direct successor}.
\end{definition}
If $\lambda : V \to \Z$ is effectively decreasing with respect to $W
\subset V$, it is important that $\lambda$ is monotone on $V$ and not
only on $W$. However it is enough to define $\lambda$ on $W$ only,
since it can be extended by $\lambda(v) \defi \min \lambda(W) - 1$ for
all $v \in V \setminus W$.

The following lemma is the main tool to bound the expected path length
of $p$.

\begin{theorem} \label{thm:ubound-a}
  Let $D = (V,A)$ be an arbitrary AUSO of a simple 4-polytope, let $\vst
  \in V$ be an arbitrary starting vertex for random edge, and let $p$ be
  the random path from $\vst$ to the global sink $\s$ defined by random
  edge starting at $\vst$.

  Let $\lambda : V \to \Z$ be an effectively decreasing function then
  $\EXP[\length(p)] \le 4 \lvert\lambda(V)\rvert$.
\end{theorem}

For simple $d$-polytopes one can prove an upper bound of $d
\#\lambda(V)$. 

\begin{proof}
  Set $V_i \defi \{\,v \in V \;:\; \lambda(v) = i \,\}$.  Then we can
  write
  \begin{equation} \label{eq:ubound:eq1}
    \EXP(p) = 
    \sum_{i \in \lambda(V)} \;
    \left(
      \sum_{e \in \incut(V_i)} p(e) \;+\; 
      \sum_{e \in E(V_i,V_i)} p(e)
    \right).
  \end{equation}
  
  If $\sum_{e \in \incut(V_i)} p(e) > 1$, there must be a (directed)
  path from a vertex in $V_i$, leaving $V_i$ and reentering $V_i$ at a
  different vertex.  But since $\lambda$ is monotone there cannot be
  such a path.  Thus
  \begin{equation} \label{eq:ubound:eq2}
    \sum_{e \in \incut(V_i)} p(e) \leq 1.
  \end{equation}
  
  In the graph of a simple 4-polytope, every vertex has degree four.
  Furthermore every vertex $v \in V_i$ has an outgoing edge $(v,w)$
  leaving $V_i$. Thus after the random edge step at $v$ the set $V_i$
  is left with probability at least $\tfrac 14$ and not revisited in
  any following random edge step.  Thus
  \begin{equation} \label{eq:ubound:eq4}
    \sum_{e \in E(V_i,V_i)} p(e) 
    \leq \sum_{j=1}^{\#V_i} \left(1-\tfrac14\right)^j
    \leq 3.
  \end{equation}
  Combining \eqref{eq:ubound:eq2} and \eqref{eq:ubound:eq4} with
  \eqref{eq:ubound:eq1} completes the proof:
  \begin{equation*}
    \EXP(p) = 
      \sum_{i \in \lambda(V)} \; \left(\,
      \sum_{e \in \incut(V_i)} p(e) \,+\,
      \sum_{e \in E(V_i,V_i)} p(e) \,\right)
      = 4\,\lvert\lambda(V)\rvert. \qedhere
    \end{equation*}
\end{proof}

Global functions $\lambda$ may be obtained from a partition
$\Pi = \{V_1, \ldots, V_\ell\}$ of the vertex set $V$ and
local functions $\lambda_1, \ldots, \lambda_\ell$.  The partition must
be ``compatible'' with the directed underlying graph in the sense that
$D/\Pi$ is acyclic. 

\begin{theorem} \label{thm:runtime}
  Let $D = (V,A)$ be an arbitrary AUSO of a simple 4-polytope, let $\vst
  \in V$ be an arbitrary starting vertex for random edge, and let $p$ be
  the random path from $\vst$ to the global sink $\s$ defined by random
  edge starting at $\vst$.

  Let $\Pi = \{V_1, \ldots, V_\ell\}$ be a partition of the vertex set
  $V$ and let $\lambda_1, \ldots, \lambda_\ell$ be effectively
  decreasing functions (with respect to $V_i$) $\lambda_i : V_i \to
  \Z$ .  Suppose that $D/\Pi$ is an acyclic digraph.  Then
  $\EXP[\length(p)] = O(\,\sum_{i=1}^\ell\#\lambda(V_i)\,)$.
\end{theorem}
\begin{proof}
  W.l.o.g. we assume that the order of the sets $V_1, \ldots, V_\ell$
  is a topological ordering of the corresponding vertices in $G/\{V_1,
  \ldots, V_\ell\}$.  Then we can use the function $\lambda : V \to
  \Z$ defined as
  \begin{equation*}
    \lambda(v) \;\defi\;
    \left( \sum_{i=1}^{\operatorname{ind}(v)-1} \max \lambda_i(V_i)
    \right) 
    \,+\, \lambda_{\operatorname{ind}(v)}(v) \,+\,
    \left( \sum_{i=\operatorname{ind}(v)+1}^\ell \min \lambda_i(V_i)
    \right).
  \end{equation*}
  where $\operatorname{ind}(v) \defi i \in [0,\ell]$ with $v \in V_i$.
  Since $G/\{V_1, \ldots, V_\ell\}$ is acyclic, the monotonicity of
  $\lambda$ follows from the monotonicity of the $\lambda_i$.
\end{proof}

Our definition of an effectively decreasing function is equivalent to
a partition
\begin{equation*}
  \Pi = \{W_1, W_2, \ldots, W_k\}
\end{equation*}
of $V$ such that $D / \Pi$ is acyclic, the numbering of the $W_i$ is a
topological ordering of the vertices in $D / \Pi$, and for each $v \in
W_i$ there is a $w \in W_j$ with $(v,w) \in A(D)$ and $i < j$.
In the light of this equivalence Theorem~\ref{thm:runtime} is just a
reformulation of Theorem~\ref{thm:ubound-a} combining monotone
functions with decreasing direct successors and partitions of the
above type. Using the combined formulation of Theorem~\ref{thm:runtime} 
is more comfortable for the proof of the main 
theorem in Section~\ref{sect:analysis}.

\subsection {Acyclic Unique Sink Orientations}
\label{sect:prelims:ausos}

From the graph $G(P)$ of a $d$-Polytope $P$ we get a directed graph
$D$ by assigning each edge of $G$ an orientation. $D$ is called a
\emph{linear orientation} if there exists a realization of $P$ in
$\R^d$ and a linear function $\phi:\R^d \to \R$ such that each edge
$\{v,w\}$ is oriented from $v$ to $w$ if and only if $\phi(v) >
\phi(w)$.
It is not known which combinatorial properties of an oriented
polytopal graph characterize linear orientations of the underlying
polytope.
                                
\emph{Acyclic unique sink orientations} are a purely combinatorial
model of orientations of polytopal graphs.  
\begin{definition} \label{def:auso} 
  Let $D$ be an orientation of a polytopal graph $G(P)$. Then $D$ is
  an \emph{acyclic unique sink orientation (AUSO)} if $D$ is acyclic
  and for every nonempty face $F \subseteq P$ the induced subgraph
  $D[F]$ has a unique sink.
\end{definition}
Every linear orientation is an AUSO, but not vice versa.  Thus AUSOs
are a more general model than linear orientations.
If $P$ is simple, then it suffices to require that only all 2-faces of
$P$ have a unique sink (see \cite{JKK02}).  There are two important
properties which follow from this fact for AUSOs on simple polytopes.
First there are also unique sources in every non-empty face of $P$,
and secondly the reverse orientation of an AUSO is again an AUSO.

\subsection {Combinatorics of Dual Cyclic Polytopes}
\label{sect:prelims:dcp}

The cyclic 4-polytope on $n$ vertices is defined as in 
\cite[p. 11]{Ziegler98}:
\begin{equation*}
  C(n) \defi \conv \left\{ \; (i, i^2, i^3, i^4)^\tr \in \R^4 : i \in
    \{0, 1, \ldots, n-1 \} \; \right\}.
\end{equation*}

All points $(i, i^2, i^3, i^4)^\tr$ are vertices of $C(n)$. We define
$C^\polar(n)$ to be the (combinatorial) polar of $C(n)$.
The combinatorics of the cyclic polytopes and thus of their duals are
given by Gale's evenness condition.

\begin{theorem}[Gale's evenness condition]
  \label {thm:gale}
  $C^\polar(n)$ is a simple polytope.  Let $f_i$ be the facet of
  $C^\polar(n)$ corresponding to the vertex $(i, i^2,i^3,i^4)$ of
  $C(n)$. Then a $4$-subset $S \subset \{0,1, \ldots, n-1\}$
  corresponds to a vertex of $C^\polar(n)$ if and only if the
  following ``evenness condition'' is satisfied:
  
  If $i < j$ are not in $S$, then the number of $k \in S$ between $i$
  and $j$ is even:

  \begin{equation}\label{eq:evenness}
    2 \;\Big|\; \#\{\: k \;:\; k \in S,\: i < k < j \:\} \quad\quad \text{for $i,j
      \notin S$}
  \end{equation}
\end{theorem}

\begin{figure}
  \begin{center}
    \begin{overpic}[tics=5,scale=.5]{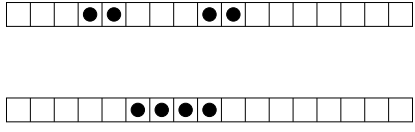}
    \end{overpic}
  \end{center}
  \caption{ \label{fig:dcp} Gale's evenness Condition: Two vertices of
    $C^\polar(17)$ illustrated by their facet incidences. The first
    one is incident to exactly two 2-faces of $\F$. The second one is
    incident to three 2-faces of $\F$.}
\end{figure}

See e.g. \cite[p. 14]{Ziegler98} for a proof.
This immediately leads to a complete description of the combinatorics
of $C^\polar(n)$.
We define the $n$ 2-faces $F_i$ to be those incident to the facets
$f_i$ and $f_{i + 1 \mod n}$. Set $\F = \{ F_0, F_1, \ldots,
F_{n-1}\}$.

Each 2-face $F \in \F$ has $(n-2)$ vertices and $C^\polar(n)$ has
$n(n-3)/2$ vertices.  A vertex is either incident to exactly two
2-faces in $\F$ or to exactly three 2-faces in $\F$ 
(see Figure~\ref{fig:dcp}).
Every vertex is uniquely determined by
\begin{align}
  \min(v) \defi &\min \{ \, i \in [0,n-1] \; : \; v \in F_i \,\} \label{eq:min}\\
  \max(v) \defi &\max \{ \, i \in [0,n-1] \; : \; v \in F_i \,\}.
  \label{eq:max}
\end{align}

We call a pair $(F_i,F_j)$ neighbors if and only if $j \equiv i + 1
\mod n$. Thus two neighboring 2-faces $(F_i,F_j)$ span the facet
$f_{i+1 \mod n}$ in the sense that every vertex of $f_{i+1}$ lies in
$F_i$ or $F_j$.
Two neighbors $(F_i,F_j)$ intersect in an edge. Thus the
facets are wedges over $(n-2)$-gons.

\begin{figure}
  \begin{center}
    \begin{overpic}[tics=5,scale=.5]{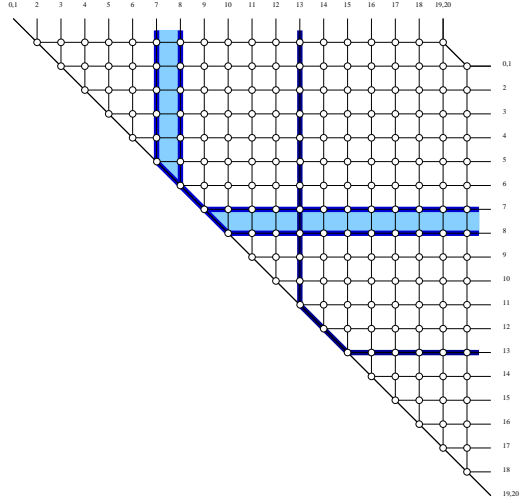}
    \end{overpic}
  \end{center}
  \caption{ \label{fig:dcp-graph}
    The graph of $C^\polar(21)$. The 2-faces $F_7$, $F_8$, and
    $F_{13}$ are indicated. The 2-faces $F_7$ and $F_8$ span a facet,
    which is also indicated.  See also \cite{HZ02} for a description
    of the graphs of dual cyclic 4-polytopes.}
\end{figure}

Furthermore we conclude from Gale's evenness condition~\ref{thm:gale}
that the facets can be renumbered in the following way.  We can choose
a facet to be the first one $f_0$ and we can reverse the numbering
keeping $f_0$, i.e. making $f_{n-1}$ the second facet and $f_1$ the
last. This corresponds to the same changes in the numbering of the
2-faces in $\F$.

\begin{definition}
  Let $F_i \in \F$ be a 2-face. We define the following vertex subsets
  of $F_i$.
  \begin{align}
    F_i^V &\defi \left\{ v \in F_i : \max(v) = i \right\} \\
    F_i^H &\defi \left\{ v \in F_i : \min(v) = i \right\}
  \end{align}
  
  A vertex $v \in V$ is called \emph{vertical} with respect to $F_i$
  if $\max(v) = i$, \emph{horizontal} with respect to $F_i$ if
  $\min(v) = i$.
  A source $\q_i$ of the 2-face $F_i$ is called \emph{vertical} if $\q_i
  \in F_i^V$, \emph{horizontal} if $\q_i \in F_i^H$,
  \emph{intermediate} otherwise.
\end{definition}
The terms vertical and horizontal refer to the vertically and
horizontally drawn parts of the 2-faces in the figures 
(e.g. Figure~\ref{fig:dcp-graph}).

\subsection {AUSOs on Dual Cyclic Polytopes}
\label{sect:prelims:dcpausos}

Let $G = G(\,C^\polar(n)\,) = (V,E)$ be the graph of a dual cyclic
4-polytope with $n$ facets and let $D = (V,A)$ be an AUSO of $G$.
Furthermore let $\q_i$ be the source and $\s_i$ the sink of the 2-face
$F_i \in \F$.
Consider two neighbors $F_i,F_j \in \F$ ($i\in \{0,1,\ldots, n-1\}$
and $j = i+1 \mod n$) and their sources $\q_i$ and $\q_j$.  $F_i$ and
$F_j$ span a 3-face $f$ of $P$,
where either $\q_i$ or $\q_j$ is the source of $f$. Thus there is a
directed path $\gamma^\q_{i,j}:\q_i \dtoorfrom \q_j$ in $D$ from $\q_i$ to
$\q_j$ or vice versa. 
We can concatenate these paths for all $i$ and obtain an (undirected)
cycle $c_\q$ which passes through all sources $\q_0, \q_1, \ldots,
\q_{n-1}$.
\begin{equation*}
  c_\q \defi \gamma^\q_{0,1} \circ \gamma^\q_{1,2} \circ \ldots \circ
  \gamma^\q_{n-2,n-1} \circ \gamma^\q_{n-1,0}
\end{equation*}

We can apply the same procedure to the sinks $\s_i$ and $\s_j$ of the
2-faces $F_i$ and $F_j$ yielding a directed path $\gamma^\s_{i,j}:\s_i
\dtoorfrom \s_j$ and thus a cycle $c_\s$ passing through the sinks of
all 2-faces in $\F$.
The next propositions state properties of the paths $\gamma^\q_{i,j}$,
$\gamma^\s_{i,j}$ and the cycles $c_\q$, $c_\s$.
The results are stated for sources only, but they can easily be
transformed to sinks by reversing the orientation of all edges. Note
that this also exchanges the functions $\preds$ and $\succs$.

\begin{proposition} \label{prop:ij-paths}  
\label{lem:ij-paths}
\label{lem:hypo-edges}
Given two neighboring 2-faces $F_i$, $F_j$, there is a directed
path $\gamma^\q_{i,j}:\q_i \dtoorfrom \q_j$ with the following
properties.\\
If $\gamma^\q_{i,j}:\q_i \dto \q_j$, then $\gamma^\q_{i,j} \cap F_j =
\{\q_j\}$
and $E(\gamma^\q_{i,j} \cap F_i, F_j) \subset \outcut(F_i)$.\\
If $\gamma^\q_{i,j}:\q_j \dto \q_j$, then $\gamma^\q_{i,j} \cap F_i =
\{\q_i\}$ and $E(\gamma^\q_{i,j} \cap F_j, F_i) \subset \outcut(F_j)$.\\
Furthermore, $\gamma^\q_{i,j}$ does not traverse the edge $F_i \cap
F_j$. 
\end{proposition}

\begin{figure}
  \begin{center}
    \begin{overpic}[tics=5,scale=.8]{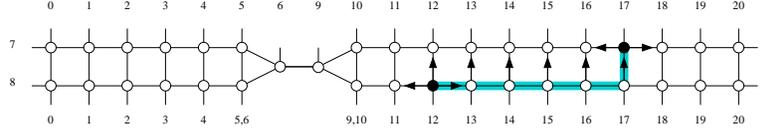}
    \end{overpic}
  \end{center}
  \caption{ \label{fig:gamma_ij} The two neighboring 2-faces $F_7$ and
    $F_8$ of $C^\polar(21)$ and the path $\gamma^\q_{7,8}$ in this case
    directed from $\q_8$ towards $\q_7$. }
\end{figure}

\begin{proof}
  See Figure~\ref{fig:gamma_ij} for an illustration of this proof.

  Let $f$ be the facet spanned by the neighboring 2-faces $F_i$ and
  $F_j$. $f$ is a simple 3-polytope and the induced subgraph $D[f]$
  has a unique source and sink.
  Let $\q_i$ and $\q_j$ be the sources of the 2-faces $F_i$ respectively
  $F_j$. Since all vertices of $f$ are vertices of $F_i$ or $F_j$
  either $\q_i$ or $\q_j$ is the unique source in $D[f]$.  Thus there
  must be a directed path from $\q_i$ or $\q_j$ to the other one.
  Assume w.l.o.g. that $\q_i$ is the source of $f$ and thus there is a
  path $\gamma :\q_i \dto \q_j$. For all those paths $\gamma \cap F_j =
  \{\q_j\}$ must hold. Otherwise there would be a directed cycle.
  Since $f$ is simple
  $\gamma$ can reach $\q_j$ only by its unique in-edge, thus only via
  its unique predecessor $v \in F_i$. And there is only one directed
  path joining $\q_i$ and $v$ without using the edge $F_i \cap F_j$.
  Thus $\gamma$ is unique and we define $\gamma^\q_{i,j} \defi \gamma$.
  
  All edges in $E(\gamma^\q_{i,j} \cap F_i, F_j)$ leave $F_i$, i.e. are
  oriented from $F_i$ to $F_j$, since an edge $e = (x,y) \in
  E(\gamma^\q_{i,j} \cap F_i, F_j)$ with $x \in F_j$ and $y \in F_i$
  would immediately imply that there is a directed cycle
  $y\stackrel{\gamma^\q_{i,j}}{\rightsquigarrow} \q_j \rightsquigarrow x
  \stackrel{e}{\to} y$.
\end{proof}

The next two propositions use the following easy observation.  Every
pair of 2-faces $(F_i,F_j)$ (not necessarily neighboring) contains at
least one vertex $v' \in F_i \cap F_j$ in its intersection. Thus there
is a directed path from $\q_i$ to $\s_j$ (via $v'$).

\begin{proposition}\label{prop:qs-paths}
\label{cor:qs-paths}
Let $F_i$ and $F_j$ be two neighboring 2-faces, then the paths
$\gamma^\q_{i,j}$ and $\gamma^\s_{i,j}$ do not intersect.
\end{proposition}
\begin{proof}
  Suppose that $\gamma^\q_{i,j}$ is directed from $\q_i$ to $\q_j$. 
  By the above observation there are paths $\omega : \q_j \dto \s_j$ and
  $\omega' : \q_j \dto \s_i$. Thus no matter how the path
  $\gamma^\s_{i,j} : \s_i \dtoorfrom \s_j$ is directed, there is a
  directed cycle if $\gamma^\q_{i,j} \cap \gamma^\s_{i,j} \not=
  \emptyset$. 
\end{proof}

\begin{proposition} \label{prop:qs}
\label{thm:qs}
The cycle $c_\q \defi \gamma^\q_{0,1} \circ \gamma^\q_{1,2} \circ \ldots
\circ \gamma^\q_{n-2,n-1} \circ \gamma^\q_{n-1,0}$ has one or two sinks
$\mqa, \mqb \in \{\q_1,\q_2,\ldots,\q_n\}$.
Not both $\mqa$ and $\mqb$ can be sinks of 2-faces in $\F$. And if
$\mqa$ is a sink of a 2-face $F_i \in \F$ then $\mqb \in \preds(\q_1)$.
\end{proposition}
\begin{proof}
  First we show, that every source of $c_\q$ must be a global source (a
  source of $D$).  The source of $c_\q$ must be a source $\q_i$, since
  all $\gamma^\q_{i,j}$ are directed.  
  
  If $\q_i$ is a source of $c_\q$, it must be the source of
  $\gamma^\q_{i-1\mod n,i}$ and $\gamma^\q_{i,i+1 \mod n}$. Thus it must
  be the source of the two 3-faces $f_i$ and $f_{i+1 \mod n}$.  With
  $f_i$ being spanned by the two 2-faces $F_{i-1 \mod n}$, $F_i$ and
  $f_{i+1 \mod n}$ being spanned by $F_i$, $F_{i+1 \mod n})$.  
  And thus $\q_i$ has at least four out edges.
  But since $C^\polar(n)$ is simple, these are all edge of $\q_i$. Thus
  $\q_i$ is the global source.
  
  Since $D$ has a unique global source $\q$, only the vertex $\q$ can be
  a source of $c_\q$. So it remains to show how many times $\q$ can be
  traversed by $c_\q$.
  If $\q$ is contained in exactly two 2-faces of $\F$, then $\q$ is
  traversed twice.
  If $\q$ is contained in exactly three 2-faces of $\F$, then these
  2-faces are of the form $F_i$, $F_{i+1 \mod n}$, and $F_{i+2 \mod
    n}$. Thus $\q$ is traversed once.
  And thus $c_\q$ has one or two sources (and of course as many sinks
  as sources).
  
  The remaining facts that not both sources can be sinks of 2-faces in
  $\F$ and that if one is such a sink, it is contained in the
  predecessors of the other, are just a simple consequence of the
  above observation, that for all pairs of 2-faces there is a directed
  path from the source to the sink of the other one.
\end{proof}

\begin{proposition} \label{prop:interm-q}
\label{lem:interm-q}
If $\q_i$ is an intermediate source, then $\q_i \in \{\mqa, \mqb, \q\}$.
\end{proposition}
\begin{figure}
  \begin{center}
    \begin{overpic}[tics=5,scale=.8]{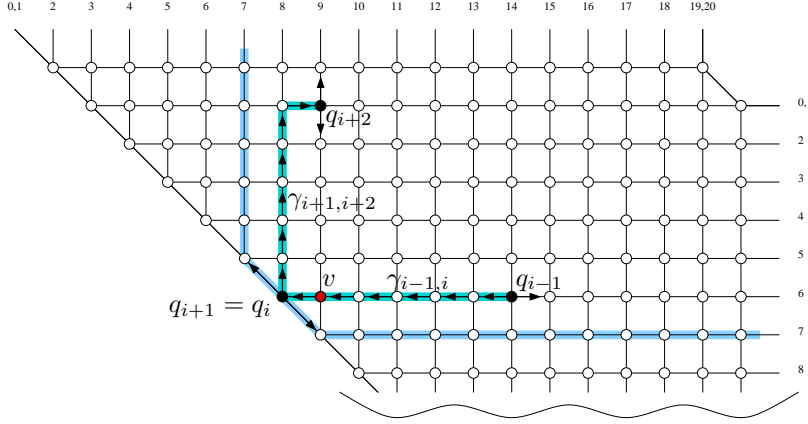}
      \put(20,13.2){$\q_{i+1} = \q_i$}
      \put(39,16.5){$v$}
      \put(63,16.5){$\q_{i-1}$}
      \put(39.4,36.8){$\q_{i+2}$}
      \put(47,16.5){$\gamma_{i-1,i}$}
      \put(34.7,26.5){$\gamma_{i+1,i+2}$}
    \end{overpic}
  \end{center}
  \caption{ \label{fig:interm-q} 
    Illustrating the proof of Proposition~\ref{lem:interm-q} in
    $C^\polar(21)$ with $i=7$.}
\end{figure}
\begin{proof}
  W.l.o.g. assume that $1 \le i \le n-3$ for the sake of not having to
  write $\mod n$ in all the following indices.  Assume that the
  intermediate source $\q_i$ is neither source nor sink of $c_\q$, i.e.
  $\q_i \notin \{\mqa, \mqb, \q\}$.
  Then we may assume w.l.o.g. (since $\q_i \neq \q$) that
  $\gamma_{i-1,i} : \q_{i-1} \dto \q_i$, i.e. $\gamma_{i-1,i}$ is
  directed from $\q_{i-1}$ to the intermediate source $\q_i$. And since
  $\q_i$ is not a sink of $c_\q$, it follows that $\q_{i+1} = \q_i$ and
  that the path $\gamma_{i+1,i+2} : \q_{i+1} \dto \q_{i+2}$ is directed
  towards $\q_{i+2}$.
  
  Now observe that $\q_i$ is the unique vertex in $F_{i+1} \cap
  F_{i-1}$, since $(i+1)-(i-1) \ge 2$.  Let $v \in F_{i-1} \cap
  F_{i+2}$.  By the same argument as before the vertex $v$ is unique,
  since $(i+2)-(i-1) \ge 2$. Thus $v$ is a neighbor of $\q_i$ on
  $F_{i-1}$.
  The path $\gamma_{i-1,i}$ must contain at least one neighbor of
  $\q_i$ on $F_{i-1}$.  Since $\gamma_{i-1,i}$ cannot traverse the edge
  $F_i \cap F_{i-1}$, $v \in \gamma_{i-1,1}$ must hold (see also
  Figure~\ref{fig:interm-q}).
  Thus we found a direct cycle $(v,\q_i) \circ \gamma_{i+1,i+2} \circ
  \q_{i+2} \rightsquigarrow v$.  A contradiction and thus $\q_i \in
  \{\mqa, \mqb, \q\}$ must hold.
\end{proof}

Note that if $\q_i = \q$ is intermediate, then $\q_{i-1} = \q_i = \q_{i+1}$
and $\q_{i-1}$ is horizontal while $\q_{i+1}$ is vertical.

\subsection {Intersecting Paths}
\label{sect:prelims:planarity}

In this section we consider a different type of graphs
to keep the results more general.

\begin{figure}
  \begin{center}
    \begin{overpic}[tics=5,scale=.75]{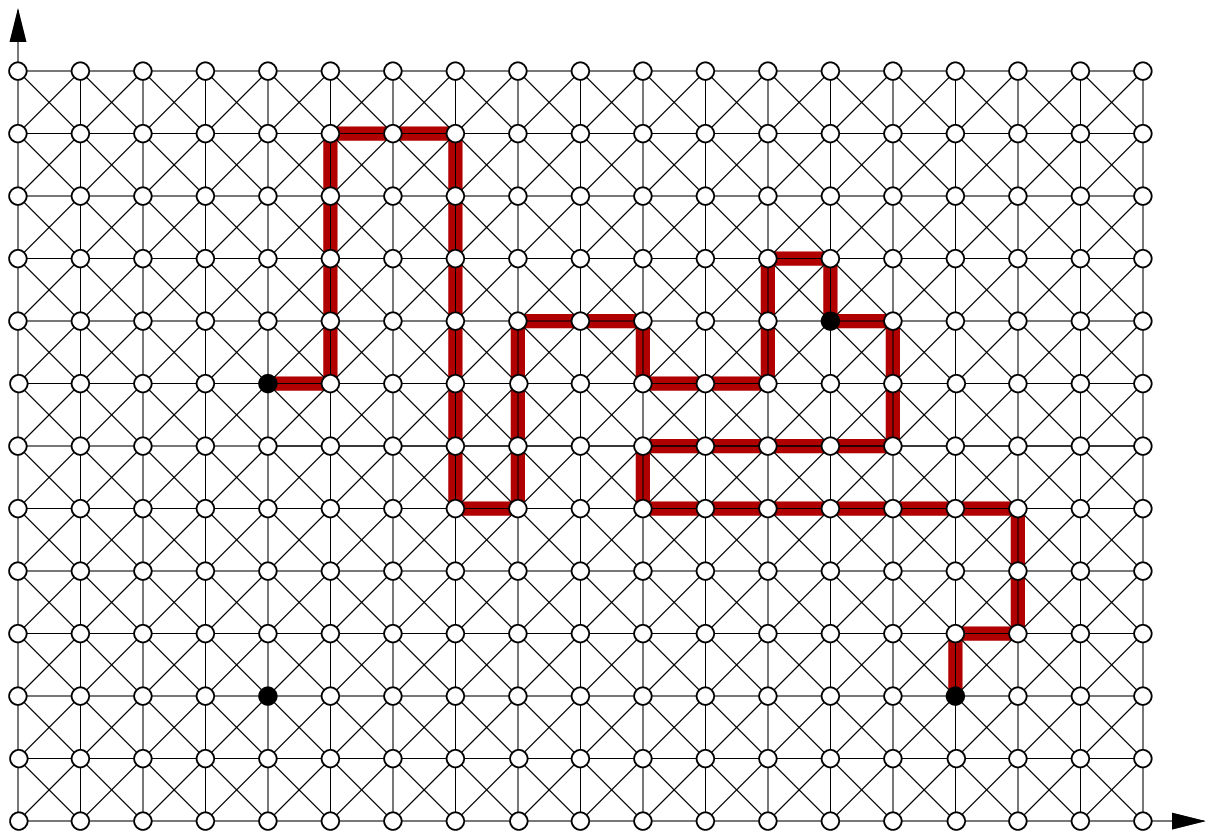}
      \put(4,1){$0$} \put(90.5,1){$n-1$} \put(102,1){$x$} \put(0,4){$0$}
      \put(-6,64){$m-1$} \put(0,68){$y$} 
      \put(24.7,1){$i$}    %
      \put(69.3,1){$r$}    %
      \put(79.2,1){$k$}    %
      \put(0.5,44.3){$s$}  %
      \put(0,39){$\ell$}   %
      \put(0,14){$j$}      %
    \end{overpic}
  \end{center}
  \caption{ \label{fig:gridgraph-b} A fence with a horizontal and
    a vertically monotone path illustrating Theorem~\ref{thm:paths}}
\end{figure}

\begin{definition}
  A graph $G = (V,E)$ is called a \emph{fence} if there are $n,m
  \in \N$ such that
  \begin {align}
    V &= \{0,1,\ldots,n-1\} \times \{0,1,\ldots,m-1\}  \label
    {eq:gridgraph:V} \\ 
    E &= \left\{\, \{(i,j)\,,\,(k,l) \,\} \;\big|\; \left|i-k\right|=1
      \text{ or } \left|j-l\right| = 1 \right\}
    \label {eq:gridgraph:E}
  \end{align}

  \begin{itemize}
  \item For a vertex $v = (i,j) \in [0,n-1] \times [0,m-1]$ we define
    $x(v) \defi i$ and $y(v) \defi j$ to be the \emph{horizontal}
    respectively \emph{vertical coordinates} of v.
  \item The sets $V_i \defi \{v \in V : x(v) = i \}$ and $H_i \defi \{
    v \in V : y(v) = i\}$ are called the \emph{vertical} respectively
    \emph{horizontal lines} of $G$.
  \item A directed fence is called \emph{sink-free} if each $V_i$
    and $H_i$ does not contain a sink.
  \end{itemize}
\end{definition}
Edges connect vertices with either only one coordinate or both
coordinates differing by one. The later edges are called \emph{skew
  edges}.  We can think of fences as graphs being embedded in
$\R^2$ with the obvious coordinates for the vertices and the edges
being straight lines (cf. Figure~\ref{fig:gridgraph-b}). Then only skew edges
cross other skew edges. All non-skew edges do not cross any other
edge.

\begin{definition}
  A path $\omega = v_1, v_2, \ldots, v_\ell$ in $G$ is called
  \emph{horizontally monotone} if for all $i \in [0,\ell-1]$
  \begin{gather} 
    x(v_{i+1}) \in \left\{x(v_i),x(v_i)-1\right\} \label{eq:monotonepaths:eq1}\\
    \left| x(v_i) - x(v_{i+1}) \right| + \left| y(v_i) - y(v_{i+1})
    \right| = 1 \label{eq:monotonepaths:eq2}\\
    (v,w) \in E(\omega \cap V_i, V_{i-1}) \Rightarrow v \in V_i \text{
      and } w \in V_{i-1} \label{eq:monotonepaths:eq3}
  \end{gather}
  holds (condition \eqref{eq:monotonepaths:eq3} only when considering
  directed fences).
  
  \emph{Vertically monotone} paths are defined analogously, with $x$
  being replaced by $y$ in \eqref{eq:monotonepaths:eq2} and $V_j$ being
  replaced by $H_j$ in \eqref{eq:monotonepaths:eq3}.
\end{definition}

\begin{definition}
  Let $D = (V,A)$ be a directed fence. A horizontally monotone path
  is called a \emph{source path} if for all $i$ with $V_i \cap \omega
  \not= \emptyset$ the path $\omega$ contains the source of the
  vertical line $V_i$.
  Analogously for $\omega$ being a vertically monotone path.
\end{definition}

 \begin{theorem} \label{thm:paths} \label{thm:paths-b} \label{cor:paths-b}
   Let $D$ be a subgraph of a directed sink-free fence.  %
   Let $v, w_1, w_2, \q \in V$ be three vertices with coordinates 
   \[
   v = (i,j), \quad w_1 = (k,j), \quad  w_2=(i,\ell), \text{ and }
   \quad \q = (r,s)
   \]
   with $k,r,s > i$ and $\ell,r,s > j$.
   Let $\omega_2 :\q \dto w_2$ be a horizontally monotone source-path
   from $\q$ to $w_2$ and $\omega_1: \q \dto w_1$ be a vertically
   monotone source-path from $\q$ to $w_1$, such that $\omega_1 \cup
   \omega_2 \subset B \defi [i,n-1] \times [j,m-1] \subset V$ and
   $\omega_1 \cap \omega_2 = \{\q\}$.

   If all skew-edges in $E(\, V \setminus B,\, (\{i, i+1, \ldots, k\} \times
   \{j\}) \cup (\{i\} \times \{j, j+1, \ldots, \ell\})\,)$ are directed
   away from $B$, i.e. they are in $\outcut(B)$, then for every $v' \in V
   \setminus B$, every directed path $\omega' : v' \dto v$ intersects
   $\omega_1 \cup \omega_2$.
 \end{theorem}

\begin{proof}
  Let $D$ be a directed sink-free fence. Consider $D$ to be
  embedded in the plane with the obvious embedding discussed above
  (which is shown in Figure~\ref{fig:gridgraph-b}). Every path or
  cycle in the graph defines a polygonal path respectively a polygonal
  curve. We will identify the path in the graph with the corresponding
  polygonal path in the embedding.
  
  Define the following (undirected) cycle $\pi$
  \[
  \pi \defi [(i,j),(i,j+1,), \ldots, (i,\ell)] \circ -\omega_2 \circ
  \omega_1 \circ [(k,j), (k-1,j), \ldots, (i,j)].
  \]
  In the embedding $\pi$ defines a simple polygonal curve.  By the
  Jordan Curve Theorem there is an interior and an exterior part of
  the plane. A point $p$ is contained in the interior, if there exists
  a ray $\rho$ (in general position) which intersects $\pi$ an odd
  number of times. In general position means, that the ray $\rho$ only
  intersects edges of $\pi$ and does not contain vertices of $\pi$.

  By condition \eqref{eq:monotonepaths:eq2} of monotone paths the
  cycle $\pi$ traverses non-skew edges only. Thus any path
  crossing $\pi$ must contain a vertex of $\pi$.

  Let $I \subset V$ be those vertices in the interior of $\pi$ or on
  $\pi$. Since $v \in I \subset B$ and $v' \notin B$ every path
  $\omega' : v' \dto v$ must intersect $\pi$ and it must traverse an
  edge in $\incut(I)$.

  Now we want to show that no edges in $\incut(I)$ are incident to
  vertices in the following two sets
  \[
  W_1 \defi \{(i,j),(i,j+1,), \ldots, (i,\ell-1)\} \quad
  W_2 \defi \{(k-1,j), (k-2,j), \ldots, (i,j)\}.
  \]
  Thus $\omega'$ must intersect $\pi$ in $\omega_1$ or $\omega_2$. 

  First, all vertices in
  \[
  \{(i+1,j),(i+1,j+1,), \ldots, (i+1,\ell)\} \,\cup\,
  \{(k,j+1), (k-1,j+1), \ldots, (i,j+1)\}
  \]
  are all either inside or on $\pi$, by the above definition of the
  interior. Thus we have to consider the edges in 
  \[
  E' \defi E(V \setminus B, W_1 \cup W_2)
  \]
  only.  We already know from the statement of the theorem that all
  skew edges in $E'$ are out-edges of $I$.  
  Thus we need to consider all edges in $E'$ along the vertical,
  respectively horizontal lines
  \[ 
  H_i, H_{i+1}, \ldots, H_{\ell-1} \quad \text{and} \quad
  V_j, V_{j+1}, \ldots, V_{k-1}.
  \]
  Let $H_c$ be one of the above horizontal lines.  By the Jordan
  Curve Theorem $H_c$ intersects $\pi$ at least twice.  Thus it must
  either intersect $\omega_1$ or $\omega_2$. If $H_c$ intersects
  $\omega_1$, then the source $\q_c$ of $H_c$ has x-coordinate $x(\q_c)
  > i$, since $\omega_1$ is a source path.  If $H_c$ intersects
  $\omega_2$ in a vertex $w$, then the edge between $w$ and the vertex
  $(x(w)-1,c)$ is directed away from $w$, by property
  \eqref{eq:monotonepaths:eq3} of the monotone path $\omega_2$.
  In both cases we conclude from $D$ being sink-free, that the
  edge between $(i,c)$ and $(i-1,c)$ is an out-edge of $I$. 
  An analogous argument holds for all vertical line $V_c$. Thus all
  Edges in $E'$ are out-edges of $I$.

  This proves the theorem for fences. It is clear that it holds
  for subgraphs of fence, too.
\end{proof}

\paragraph{Fences and graphs of dual cyclic 4-polytopes}
Now let $D = (V,A)$ or $G = (V,E)$ be the directed or undirected graph
of a dual cyclic polytope $C^\polar(n)$. 

For a map $\Phi : V \to [0,n-1]^2$ we extend $\Phi$ 
in the following obvious ways
\begin{gather*}
  \Phi(V) \defi \left\{ \Phi(v) : v \in V \right\} \\
  \begin{aligned}
    \Phi\left(\,\{v,w\}\,\right) &\defi
    \left\{\,\Phi(v),\Phi(w)\,\right\} &
    \Phi\left(\,(v,w)\,\right) &\defi \left(\,\Phi(v),\Phi(w)\,\right)\\ 
    \Phi(E) &\defi \left\{ \Phi(e) : e \in E \right\} &
    \Phi(A) &\defi \left\{ \Phi(e) : e \in A \right\} \\
    \Phi(D) &\defi \left(\,\Phi(V),\Phi(A)\,\right) &
    \Phi(G) &\defi \left(\,\Phi(V),\Phi(E)\,\right)
  \end{aligned}
\end{gather*}

\begin{definition}
  For a graph $G = (V,E)$ respectively $D = (V,A)$ an injective map
  $\Phi : V \to [0,n-1]^2$ is called a  \emph{fence
    embedding} if $\Phi(G)$ respectively $\Phi(D)$ is the subgraph of
  a fence.
\end{definition}

The following four maps define fence embeddings
of $G(\,C^\polar(n)\setminus F_1\,)$.
\begin{align*} %
  \Phi_1(v) \defi \left(\, \min(v),\max(v) \, \right) &&
  \Phi_2(v) \defi \left(\, \max(v),n-\min(v)+1 \, \right) 
\end{align*}

\section {Runtime analysis of random edge on $C^\polar(n)$}
\label{sect:analysis}

Let $n \geq 5$ be an arbitrary integer. And let $D = (V,A)$ be an
arbitrary AUSO of the graph $G(\,C^\polar(n)\,)$ of the
dual cyclic 4-polytope with $n$ facets. Let $\vst \in V$ be an
arbitrary starting vertex for random edge and let $(\,p(e)\,)_{e \in
  A}$ 
be the corresponding random path. Then we want to bound the running
time of random edge with $O(n)$, i.e. we want to show that $\sum_{e
  \in A}p(e) = O(n)$.

The line of arguments should be pretty clear by now. We will define a
constant size partition $\Pi$ of the vertex set $V$ with effectively
decreasing functions (with respect to $W$) $\lambda_W : W \to \Z$ for
each $W \in \Pi$ in order to apply Theorem \ref{thm:runtime}.  Then we
need to show that $\#\lambda_W(W) = O(n)$ for each $W \in \Pi$ to get
a linear upper bound.

We split the vertex set $V$ into three sets which will be refined later.
\begin{align}
  \Qset \defi& \preds(\mqa) \cup \preds(\mqb)\\
  \Sset \defi& \succs(\msa) \cup \succs(\msb)\\
  R \defi& V \setminus (\Qset \cup \Sset)
\end{align}

\subsection {\boldmath Investigating the Vertex Sets $\Qset$ and $\Sset$}
\label{sect:analysis:qs}

In this section we fix the 2-face $F_0 \in \F$ such that $\mqa = \q_0$
is the source of $F_0$.  This leaves us the freedom to choose between
the two possible numberings of the $F_i \in \F$ with $F_0$ fixed 
and the property that all pairs $F_i,F_{i+1}$ span a 3-face.

We would like to save some work and exploit the symmetry, 
exchanging the sets $\Qset$ and $\Sset$ by inverting the AUSO.
Inverting the AUSO also changes 
a monotone decreasing function
$\lambda_\Qset : \Qset \to \Z$ into a monotone increasing function
$\lambda_\Sset : \Sset \to \Z$ and a decreasing direct
successor becomes a decreasing direct predecessor.  Thus in order to
exploit the above symmetry we need monotone functions with decreasing
direct successors and increasing direct predecessors.  This is
formalized by conditions \eqref{eq:qs:q} respectively \eqref{eq:qs:s}.
\begin{subequations} \label{eq:qs}
\begin{align}
  \forall v \in \Qset: \; \exists (v,w) \in A: \; &
  \lambda_\Qset(v) > \lambda_\Qset(w) \mbox{ or } w \notin \Qset \label{eq:qs:q}\\
  \forall v \in \Qset: \; \exists (w',v) \in A: \; & \lambda_\Qset(w') >
  \lambda_\Qset(v) \mbox{ or } w' = q \label{eq:qs:s}
\end{align}
\end{subequations}
Thus proving equations \eqref{eq:qs} for $\lambda_\Qset$ on $\Qset$ suffices
to deal with $\Qset$ and $\Sset$, as \eqref{eq:qs:q} is the ``decreasing
direct successor'' condition for $\lambda_\Qset$ to apply Theorem
\ref{thm:runtime} on $\Qset$ and \eqref{eq:qs:s} implies the ``decreasing
direct successor'' condition to apply Theorem \ref{thm:runtime} for
$(-\lambda_\Sset)$ on $\Sset$.
Referring to \eqref{eq:qs} will always be a shortcut to refer to 
\eqref{eq:qs:q} and \eqref{eq:qs:s}.

\paragraph{Some definitions and basic properties}

Now we will settle some definitions for this section.  We assume
w.o.l.g. that $\mqa \notin \preds(\mqb)$ holds.  And if $c_\q$ has only
one sink, we will define $\mqb \defi \q$ for the following definitions
to be well defined in any case.
Define $k,\ell \in [0,n-1]$ as the coordinates of the global source
$\q$, i.e. $k \defi \min(\q)$ and $\ell \defi \max(\q)$. And define $m
\in [k,l]$ by $\q_m = \mqb$.
The cycle $c_\q$ can be split into four directed paths.  $\gamma_1$
leading from $\q$ to $\mqa$ via $\q_k,\q_{k-1}, \ldots, \q_1,\q_0$.
$\gamma_2$ leading from $q$ to $\mqa$ via $\q_\ell,\q_{\ell+1}, \ldots,
\q_{n-1},\q_0$.  $\gamma_3$ leading from $q$ to $\mqb$ via $\q_k,
\q_{k+1}, \ldots, \q_{m-1}, \q_m$.  And finally $\gamma_4$ leading from
$\q$ to $\mqb$ via $\q_\ell, \q_{\ell-1}, \ldots, \q_{m+1}, \q_m$.
Define the following ten vertex sets. (Note that $\q_{\min(v)}
\in \gamma_\alpha$ and $\q_{\max(v)} \in \gamma_\beta$ implies $\alpha
\leq \beta$.)
\begin{equation*}
  V_{\alpha,\beta} \defi \left\{\, v \in V \;:\; \q_{\min(v)} \in
    \gamma_\alpha \text{ and } \q_{\max(v)} \in \gamma_\beta \,\right\}
\end{equation*}
Set $s,t \in [0,n-1]$ to be indices, such that $\q_s \in \gamma_1 \cap
\succs(\mqb)$ and $\q_t \in \gamma_2 \cap \succs(\mqb)$
are the first such sources on $\gamma_1$ respectively $\gamma_2$. %
\label{def:qs} \label{def:qt} %

\begin{figure}
  \begin{center}
    \begin{overpic}[tics=5,scale=.75]{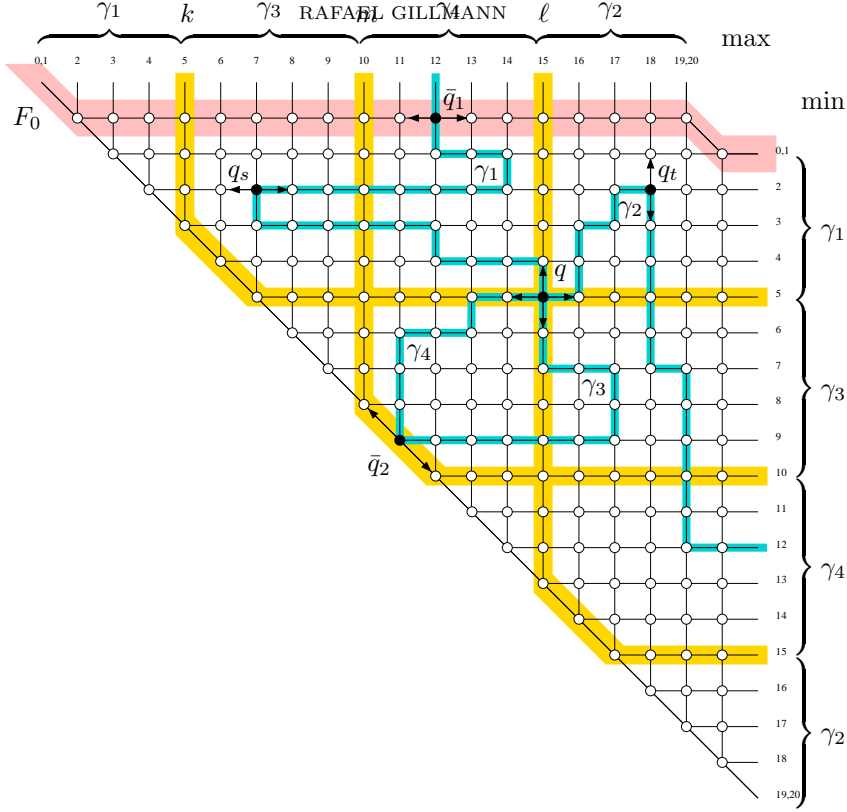}
      \put(2.5,95){\begin{tabular}[b]{c} %
          $\gamma_1$\\$\overbrace{\hspace{1.9cm}}$ %
        \end{tabular}}
      \put(20.5,95){\begin{tabular}[b]{c} %
          $\gamma_3$\\ $\overbrace{\hspace{2,35cm}}$ %
        \end{tabular}}
      \put(43,95){\begin{tabular}[b]{c} %
          $\gamma_4$\\$\overbrace{\hspace{2,35cm}}$ %
        \end{tabular}}
      \put(65,95){\begin{tabular}[b]{c} %
          $\gamma_2$\\$\overbrace{\hspace{2cm}}$ %
        \end{tabular}}

      \put(97.5,81){\begin{tabular}[t]{lr} %
          \begin{rotate}{-90}$\overbrace{\hspace{1.9cm}}$\end{rotate} %
           & \raisebox{-1.05cm}{$\gamma_1$} \end{tabular}%
       }
      \put(97.5,63){\begin{tabular}[t]{lr} %
          \begin{rotate}{-90}$\overbrace{\hspace{2.35cm}}$\end{rotate} %
           & \raisebox{-1.275cm}{$\gamma_3$} \end{tabular}%
       }
      \put(97.5,40.5){\begin{tabular}[t]{lr} %
          \begin{rotate}{-90}$\overbrace{\hspace{2.35cm}}$\end{rotate} %
           & \raisebox{-1.275cm}{$\gamma_4$} \end{tabular}%
       }
      \put(97.5,18){\begin{tabular}[t]{lr} %
          \begin{rotate}{-90}$\overbrace{\hspace{2cm}}$\end{rotate} %
           & \raisebox{-1.1cm}{$\gamma_2$} \end{tabular}%
       }
       
       \put(22,98){$k$} \put(44,98){$m$} \put(67,98){$\ell$}
       
       \put(45.5,41.5){$\mqb$} \put(55,87.5){$\mqa$} \put(69,66){$q$}
       \put(28,78.5){$\q_s$} \put(82,78.5){$\q_t$}
       
       \put(59,78.5){$\gamma_1$} \put(77,74){$\gamma_2$}
       \put(72.5,51.5){$\gamma_3$} \put(50.5,56){$\gamma_4$}
       
       \put(1,85){$F_0$}
       
       \put(90,95){$\max$} \put(100,87){$\min$}
     \end{overpic}
  \end{center}
  \caption{ \label{fig:dcp-aufteilung} The graph of $C^\polar(21)$
    with the sets $V_{i,j}$ and the paths $\gamma_1, \gamma_2,
    \gamma_3, \gamma_4$.}
\end{figure}

Before we continue defining a suitable partition of $\Qset$, we will
prove some basic facts about the geometric setting. (c.f.
Figure~\ref{fig:dcp-aufteilung}).

\begin{lemma} \label{lem:monotone-paths}
  If $\q_i \in \preds(\mqa)$ is the source of the 2-face $F_i$, $\q_i
  \in \gamma_\alpha$, and $\q_i \notin \{\mqa, \mqb\}$, then
  \begin{align*}
    \min(\q_i) &= i \quad \text{ for $\alpha = 1,3$} &&
    \text{i.e. $\q_i$ is \emph{horizontal}} \\
    \max(\q_i) &= i \quad \text{ for $\alpha = 2,4$} && \text{i.e.
      $\q_i$ is \emph{vertical}}.
  \end{align*}
\end{lemma}

\begin{proof}
  If there are indices $i,j$ such that $F_i$ and $F_j$ are neighbors
  and $\q_i,\q_j \in \preds(\mqa)$ and $\min(\q_i) = i$ and $\max(\q_j) =
  j$, then the path $\gamma_{i,j}$ must either intersect $F_0$ or
  contain the edge $e = F_i \cap F_j$.  But either one is impossible,
  since it would induce a cycle.  By Lemma \ref{lem:hypo-edges}
  respectively $\gamma_\alpha\vert_{[q,\q_i]} \cap F_0 = \emptyset$
  follows from $\q_i \in \preds(\mqa)$.
  Thus $\gamma_\alpha\vert_{[q,\q_i]}$ cannot contain vertical and
  horizontal sources. And since only $\mqa$ and $\mqb$ can be
  intermediate sources by Lemma \ref{lem:interm-q} and $\q_i \notin
  \{\mqa,\mqb\}$ the path $\gamma_\alpha\vert_{[q,\q_i]}$ contains
  either vertical or horizontal sources only.
  The paths $\gamma_1$ and $\gamma_3$ start at a horizontal source.
  While the paths $\gamma_2$ and $\gamma_4$ start at a vertical
  source.
\end{proof}

\begin{lemma} \label{lem:mqb}
  If $v \in \preds(\mqa)$ and $ v \in V_{11} \cup V_{13} \cup V_{33}
  \cup V_{22} \cup V_{24} \cup V_{44}$
  then $v \in \succs(\mqb)$ and thus $\mqb \in \preds(\mqa)$.
\end{lemma}

\begin{proof}
  Set $i \defi \min(v)$. Suppose $v \in V_{11} \cup V_{13} \cup
  V_{33}$ thus $\q_i \in \gamma_\alpha$ with $\alpha = 1$ or $\alpha =
  3$.  There is a path $\gamma' : \q_i \dto v$ with $\gamma' \subset
  F_i$.
  
  If $\max(\q_i) > m$, then either $\gamma' \cap F_0 \neq \emptyset$ or
  $\gamma' \cap F_m \neq \emptyset$. While the first is impossible
  since it imposes a cycle in $D$, the latter proves $v \in
  \succs(\mqb)$.
  
  If $\max(\q_i) < m$ then set $\q_j \in \gamma_\alpha$ to be the first
  source on $\gamma_\alpha$ with $\max(\q_j) < m$.  Set $\hat{\jmath}
  \defi j+1$ if $\alpha = 1$ and $\hat{\jmath} \defi j-1$ if $\alpha =
  3$, i.e. $\q_{\hat{\jmath}}$ is the source before $\q_j$ on
  $\gamma_\alpha$.  Since $q=\q_k \in \gamma_1 \cap \gamma_3$ is the
  first vertex of $\gamma_1$ and $\gamma_3$ and $\max(\q_k) = l > m$,
  $\q_j \not= q$ holds. And thus $\q_{\hat{\jmath}}$ exists and
  $\max(\q_{\hat{\jmath}}) \ge m$.  Thus $\gamma_{\hat{\jmath},j} \cap
  F_m \neq \emptyset$, since $\gamma_{\hat{\jmath},j} \cap F_0 =
  \emptyset$. And since $\gamma_{\hat{\jmath},j} \subset
  \gamma_\alpha\vert_{[q,\q_i]}$ it follows that $v \in \succs(\mqb)$.
  (Note that if $\alpha = 3$ the condition $\gamma_{\hat{\jmath},j}
  \cap F_m \neq \emptyset$ imposes a cycle in $D$ and thus is a
  contradiction.)

  The case $v \in V_{22} \cup V_{24} \cup V_{44}$ is proven by
  symmetry. Reversing the numbering of the $F_i \in \F$ while keeping
  $F_0$ fixed exchanges the path $\gamma_1$ and $\gamma_2$
  respectively $\gamma_3$ and $\gamma_4$.
\end{proof}

\paragraph{Defining a suitable partition}

\begin{figure}
  \begin{center}
    \begin{overpic}[tics=5,scale=.5]{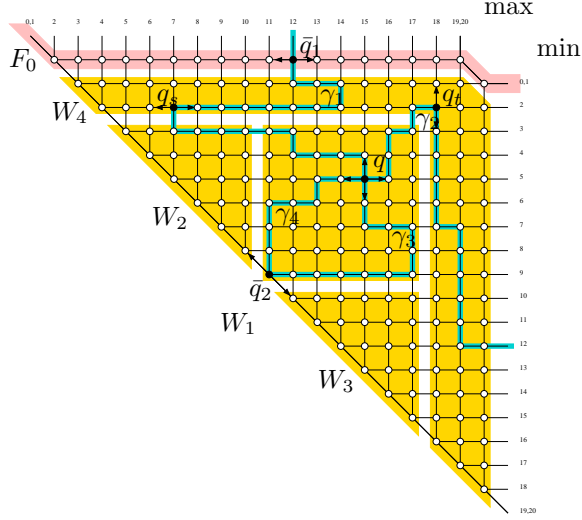}
       \put(45.5,41.5){$\mqb$} \put(55,87.5){$\mqa$} \put(69,66){$q$}
       \put(28,78.5){$\q_s$} \put(82,78.5){$\q_t$}
       
       \put(59,78.5){$\gamma_1$} \put(77,74){$\gamma_2$}
       \put(72.5,51.5){$\gamma_3$} \put(50.5,56){$\gamma_4$}
       
       \put(0,85){$F_0$}
       
       \put(90,95){$\max$} \put(100,87){$\min$}

       \put(8,75){$W_4$}
       \put(27,55){$W_2$}
       \put(40,35){$W_1$}
       \put(58.5,24){$W_3$}
     \end{overpic}
  \end{center}
  \caption{ \label{fig:dcp-w-aufteilung} The graph of $C^\polar(21)$
    with the sets $W_1$, $W_2$, $W_3$, and $W_4$}
\end{figure}

In view of Theorem~\ref{thm:runtime} we would like to define a
partition $\Pi_\Qset$ of $\Qset$ such that $D[\Qset]/\Pi_\Qset$ is acyclic.  To get
this partition and a useful characterization, we will first define the
following partition of the vertex-set $V \setminus F_0$ (c.f.
Figure~\ref{fig:dcp-w-aufteilung}).
\begin{align}
  W_1 &\defi \left( F_{m+1} \cup F_{m+2} \cup \ldots \cup F_{t-1}
  \right) \;\cap\; \left( F_{s+1} \cup F_{s+2} \cup \ldots \cup
    F_{m-1} \right) \\
  W_2 &\defi \left(F^V_{s+1} \cup F^V_{s+2} \cup \ldots \cup F^V_{m}
  \right) \;\cap\; \left(F^H_{s+1} \cup F^H_{s+2} \cup \ldots \cup
    F^H_m \right) \\
  W_3 &\defi \left( F^V_m \cup F^V_{m+1} \cup \ldots \cup F^V_{t-1}
  \right) \;\cap\; \left(F^H_m \cup F^H_{m+1} \cup \ldots \cup
    F^H_{t-1} \right) \\
  W_4 &\defi \left( F_1 \cup F_2 \cup \ldots \cup F_s\right) \;\cup\;
  \left(F_t \cup F_{t+1} \cup \ldots \cup F_{n-1}\right)
\end{align}
Now define the following sets.
\begin{align}
  \Qset_1 &\defi \preds(\mqa) \setminus \succs(\mqb) \\
  \Qset_2 &\defi \preds(\mqa) \cap \succs(\mqb) \setminus
  \left(\succs[\q_s] \cap \succs[\q_t]\right) \cap %
  \left( V_{11} \cup V_{13} \cup V_{33}\right)
  \\
  \Qset_3 &\defi \preds(\mqa) \cap \succs(\mqb) \setminus
  \left(\succs[\q_s] \cap \succs[\q_t]\right) \cap %
  \left( V_{22} \cup V_{24} \cup V_{44}\right)
  \\
  \Qset_4 &\defi \preds(\mqa) \cap \succs(\mqb) \cap (\succs[\q_s] \cup
  \succs[\q_t])
\end{align}
If $c_\q$ has only one sink and $\mqb = q$, then all sets $\Qset_i$ are
still well defined. And note that in this case $\q_s = \q_t = q$ holds
and thus $\Qset_1 = \{q\}$ and $\Qset_2 = \Qset_3 = \emptyset$.
If $\mqb \notin \preds(\mqa)$ then $\Qset_2 = \Qset_3 = \Qset_4 = \emptyset$ and
$\Qset_1 = \preds(\mqa)$ hold.

\begin{lemma} \label{lem:Q-partition}
  For all $i \in \{1,2,3,4\}$, $W_i \cap \preds(\mqa) = \Qset_i$
  holds. Furthermore 
  \begin{equation*}
    \Pi_\Qset \defi \{\Qset_1, \Qset_2, \Qset_3, \Qset_4\}
  \end{equation*}
  is a partition of $\Qset$ and $D[q]/\Pi_\Qset$ is acyclic.
\end{lemma}
\begin{proof}
  First note that the sets $W_1, W_2, W_3, W_4$ are a partition of $V
  \setminus F_0$ (c.f. Figure~\ref{fig:dcp-w-aufteilung}). And thus
  the sets $W_i \cap \preds(\mqa)$ are a partition of $\Qset =
  \preds(\mqa)$.

  Now we prove $W_1 \cap \preds(\mqa) = \Qset_1$.
  By the definition of $s$ and $t$ $W_4 \subset \succs(\mqb)$
  holds. By Lemma~\ref{lem:mqb} $V_{11} \cup V_{13} \cup V_{33} \cup
  V_{22} \cup V_{24} \cup V_{44} \subset \succs[\mqb]$ holds. Thus
  \begin {equation} \label{eq:Q-partition:Q1-1}
    \preds(\mqa) \setminus \succs(\mqb) \subset \left(\left( \left( V_{14}
        \cup V_{12} \cup V_{23} \cup V_{34} \right) \setminus F_m \right)
    \cup \{\mqb\} \right) \setminus (W_2 \cup W_3) = W_1
  \end{equation}
  holds.
  It suffices to prove $\incut(W_1 \cap \preds(\mqa)) = \emptyset$
  and $\outcut(\mqb) \subset \outcut(W_1)$. Then $W_1 \cap
  \preds(\mqa) \cap \succs(\mqb) = \emptyset$ and thus $W_1 \cap
  \preds(\mqa) \subset \preds(\mqa) \setminus \succs(\mqb)$ holds. 
  Suppose there is an in-edge $e = (v,w) \in \incut(W_1 \cap
  \preds(\mqa))$.  Since $\incut(\preds(\mqa)) = \emptyset$, $e \in
  \incut(W_1)$ holds.
  Set $i \defi \min (w)$ and $j  \defi \max(w)$. 
  And thus either $e \in F_i$ or $e \in F_j$.
  By the definition of $W_1$, $i \in \{ s+1, s+2,
  \ldots, m-1\}$ and $j \in \{ m+1, m+2, \ldots, t-1 \}$ holds.
  By the choice of $s$ and $t$ and \eqref{eq:Q-partition:Q1-1}
  \begin{equation*} \label{eq:Q-partition:Q1-2}
    \gamma_1\vert_{[k,s-1]},\; 
    \gamma_2\vert_{[k,t-1]},\;
    \gamma_3\cap\preds(\mqa),\;
    \gamma_4 \cap \preds(\mqa) \; \subset \;\preds(\mqa)
    \setminus \succs(\mqb)
  \end{equation*}
  and thus
  \begin{equation*}
    \q_i, \q_j \in W_1 \cap \preds(\mqa)
  \end{equation*}
  holds. Since
  \begin{equation*}
    \preds(\mqa)
    \setminus \succs(\mqb)
    \subset W_1 \cap \preds(\mqa).
  \end{equation*}

  Assume $e \in F_i$. (The argument for $e \in F_j$ is the
  same.) 
  The sink $\s_i$ of $F_i$ must be located on $F_i$ within $W_1$
  between $w$ and $\q_i$. Thus the (unique) path $\gamma':\q_i \dto w
  \subset F_i$ must pass $F_0$. But then there would be a cycle again
  (since $\q_i, w \in \preds(\mqa)$).
  If $\mqb \not= q$, $\mqb$ has only two direct successors on $F_m$
  and thus outside of $W_1$. If $\mqb = q$ then $W_1 =
  \{q\}$. Therefor $\outcut(\mqb) \subset \outcut(W_1)$ holds.

  Now we prove $W_4 \cap \preds(\mqa) = \Qset_4$.
  $W_4 \subset \succs[\q_s] \cup \succs[\q_t]$ by the definition of
  $W_4$ and thus $W_4 \cap \preds(\mqa) \subset \preds(\mqa) \cap
  (\succs[\q_s] \cup \succs[\q_t])$. 
  We will show that $\outcut(W_4 \cap \preds(\mqa)) \subset
  \outcut(\preds\mqa)$. That suffices, since $\q_s, \q_t \in W_4 \cap
  \preds(\mqa)$ implies that $(\succs[\q_s] \cup \succs[\q_t]) \cap
  \preds(\mqa) \subset W_4\cap\preds(\mqa)$.
  Suppose there is an edge $(v,w) \in \outcut(W_4 \cap \preds(\mqa)
  ) \setminus \outcut(\preds(\mqa)) $, then $v \in W_4 \cap
  \preds(\mqa)$ and $w \in (W_2 \cup W_3) \cap \preds(\mqa)$.
  (Remember that $W_1 \cap \preds(\mqa)$ has only out-edges.)  We
  assume that $w \in W_2$ since the argument for $w \in W_3$ is the
  same. If $w \in W_2$ then $v \in F_s$ and $w \in F_{s+1}$ holds by
  the construction of $W_4$ and $W_2$. And thus $\min(v) = s$,
  $\min(w) = s+1$. Choose $j$ such that $e \in F_j$, then $s < j \le
  m$ and $\q_j \in \gamma_1$ or $\q_j \in \gamma_3$. In either case $\q_j
  \in \preds(\mqa)$ and $\q_j \notin \succs(\mqb)$ holds. Thus $\q_j \in
  W_1 \cap \preds(\mqa)$ and thus $\s_j \in F_j \cap (F_{s+1} \cup
  F_{s+2} \cup \cdots \cup F_{\max(\q_j)})$ and finally %
  there is a
  directed path $\gamma' : \q_j \dto v$, $\gamma' \subset F_j$ which intersects
  $F_0$. This is again a contradiction.

  From the fact that $\Qset_2, \Qset_3 \not\subset \Qset_1 \cup
  \Qset_4$, it follows that $\Qset_i = W_i \cap \preds(\mqa)$.

  Thus $\Pi_\Qset$ is a partition of $\Qset$. Since $\incut(\Qset_1) =
  \outcut(\Qset_4) = \emptyset$ any cycle in $D[\Qset]/\Pi_\Qset$ can contain only
  the vertices $\Qset_2$ and $\Qset_3$. But since for all $v \in W_2$ and $w
  \in W_3$ both $\min(w) - \min(v) > 1$ and $\max(w) - \max(v) > 1$ hold,
  there cannot be any edges in $D$ between $W_2$ and $W_3$,
  i.e. $E(\Qset_1,\Qset_2) = \emptyset$. Thus $D[\Qset]/\Pi_\Qset$ is acyclic.
\end{proof}

\paragraph{Defining monotone functions}

Now we have a partition 
\begin{equation*}
  \Pi_\Qset = \{\Qset_1, \Qset_2, \Qset_3,
\Qset_4 \}
\end{equation*}
of $\Qset$. This is the first ingredient to
Theorem~\ref{thm:runtime}. The second ingredient are effectively
decreasing functions. As discussed at the beginning of this section we
would like to exploit symmetries in the definition of $\Qset$ and
$\Sset$, so we are looking for monotone decreasing function
$\lambda_1, \lambda_2, \lambda_3, \lambda_4$ which satisfy
\eqref{eq:qs} (w.r.t. $\Qset_i$).

Define the function $\lambda_\Qset : \Qset \to \Z$ as
\begin{equation*}
  \lambda_\Qset(v) \defi - \#\left\{ \, \preds[v] \cap \{\q_0, \q_1, \ldots,
    \q_{n-1}\} \,\right\}
\end{equation*}
the number of sources $\q_i$ in the predecessors of $v$.  This is a
monotone decreasing function on $\Qset$.
Then we set 
\begin{align}
  \lambda_1 &\defi \lambda_\Qset & \lambda_4 &\defi \lambda_\Qset.
\end{align}

Unfortunately the sets $\Qset_2$ and $\Qset_3$ may be large sets on which
\eqref{eq:qs} does not hold for $\lambda_\Qset$. Thus we need different
functions on those two sets.
We define
\begin{align}
  \lambda_2 &\defi \min && \lambda_3 \defi -\max.
\end{align}
Where $\min$ and $\max$ are the coordinate-functions.

\paragraph{\boldmath Analyzing $\lambda_2$ on $\Qset_2$ and $\lambda_3$ on
  $\Qset_3$}

It is easier to prove \eqref{eq:qs} for the functions $\lambda_2$ and
$\lambda_3$. Since it is not clear, that they both are monotone
decreasing, we have to prove the following lemma.

\begin{lemma} \label{lem:main-Q2}
  $\min:V \rightarrow \Z$ is monotone decreasing and satisfies
  conditions \eqref{eq:qs} on $\Qset_2$.
\end{lemma}

\begin{proof}
  By Lemma~\ref{lem:mqb} $\mqb \in \preds(\mqa)$. Thus applying
  Lemma~\ref{lem:monotone-paths} yields that all sources $\q_i$ in the
  paths $\gamma_1$ and $\gamma_3$ are horizontal sources.
  For a $v \in \Qset_2$ set $j \defi \max(v)$. Then $\q_j \in \gamma_1 \cup
  \gamma_3$ and there exists a $\gamma': \q_j \dto v$ with $\gamma'
  \subset F_j$.  Since $\q_j$ is a horizontal source either there is a
  $v' \in \gamma'$ with $\min(v') = j-1$ and $\max(v') = j+1$ or
  $\gamma' \cap F_0 \neq \emptyset$. Since the latter is impossible it
  follows that the first one holds. And since $v \neq \s_j$ $(w,v) \in
  A$ and $(v,w') \in A$ with $w,w' \in F_j$ and $\min(w) > \min(v) >
  \min(w')$.
\end{proof}

\begin{corollary} \label{cor:main-Q3}
  $-\max:V \rightarrow \Z$ is monotone decreasing and satisfies
  conditions \eqref{eq:qs} on $\Qset_3$.
\end{corollary}
\begin{proof}
  Reversing the numbering of the 2-faces $\F$ while keeping $F_0$
  fixed exchanges $\gamma_1$ with $\gamma_2$ and $\gamma_3$ with
  $\gamma_4$ and thus exchanges the sets $\Qset_2$ and $\Qset_3$.
\end{proof}

\paragraph{\boldmath Analyzing $\lambda_\Qset$ on $\Qset_1$ and $\Qset_4$}

It is clear that $\lambda_\Qset$ is monotone decreasing. Thus we only have
to prove, that \eqref{eq:qs} is satisfied.
This involves more technical details than proving
Lemma~\ref{lem:main-Q2}. 
We need to argue that certain paths, if they exist, must
intersect. This is done by applying Corollary~\ref{cor:paths-b} from
Section~\ref{sect:prelims:planarity}. But in order to apply
Corollary~\ref{cor:paths-b} we have to define suitable fence
embeddings first. We start by giving several fence embeddings
and proving some properties of them.

\begin{lemma} \label{lem:embeddings}
For $\alpha \in \{1,2,3,4\}$ let $\Phi$ be one of the following two fence 
embeddings depending on $\alpha$ and leaving a choice in the $y$-Coordinate.
  \begin{align}
    \alpha = 1:& \quad\left|\quad\begin{aligned}
        x(\,\Phi(w)\,) &\defi n - \min(w) + 1\\
        y(\,\Phi(w)\,) &\defi \begin{cases}
          \max(w)&\text{or}\\
          n - \max(w) + 1
        \end{cases}
      \end{aligned}\right.\\
    \alpha=2:& \quad\left|\quad\begin{aligned}
        x(\,\Phi(w)\,) &\defi \max(w)\\
        y(\,\Phi(w)\,) &\defi \begin{cases}
          \min(w)&\text{or}\\
          n - \min(w) + 1
        \end{cases}
      \end{aligned}\right.\\
    \alpha=3:& \quad\left|\quad\begin{aligned}
        x(\,\Phi(w)\,) &\defi \min(w)\\
        y(\,\Phi(w)\,) &\defi \begin{cases}
          \max(w)&\text{or}\\
          n - \max(w) + 1
        \end{cases}
      \end{aligned}\right.\\
    \alpha=4:& \quad\left|\quad\begin{aligned}
        x(\,\Phi(w)\,) &\defi n-\max(w)+1\\
        y(\,\Phi(w)\,) &\defi \begin{cases}
          \min(w)&\text{or}\\
          n - \min(w) + 1
        \end{cases}
      \end{aligned}\right.
  \end{align}

  Then $\Phi(\gamma_\alpha)$ is a horizontal monotone source-path and
  a vertical monotone source-path if $x(\,\Phi(w)\,)$ and
  $y(\,\Phi(w)\,)$ are exchanged.
\end{lemma}

\begin{proof}
  Check that each possible $\Phi$ is a fence embedding which for
  $\alpha = 2,4$ embeds the $F^H_i$ into the horizontal lines and the
  $F^V_i$ into the vertical lines of the fence and for
  $\alpha=1,3$ vice versa.
  
  From Lemma \ref{lem:monotone-paths} we can conclude, that
  $\gamma_\alpha$ has vertical or horizontal sources only.  Now check
  that the coordinates of $\Phi$ are chosen appropriately to make the
  path $\gamma_\alpha$ horizontal monotone. All properties that a
  directed monotone path needs are an immediate translation of Lemma
  \ref{lem:ij-paths} into the terms of Section
  \ref{sect:prelims:planarity}.
\end{proof}

\begin{lemma} \label{lem:embeddings-sinkfree}
  For any embedding $\Phi$ in Lemma~\ref{lem:embeddings} the graph
  $\Phi(\,D[\preds(\mqa)]\,)$ is the subgraph of a sink-free
  fence.
\end{lemma}

\begin{proof}
  As $\preds(\mqa)$ does not contain any sinks $\s_i$ of the 2-faces
  $F_i$ and as only one $F_i$ is mapped to each horizontal and
  vertical line it suffices to show that each $F_i \cap \preds(\mqa)$
  is an interval (i.e. connected).
  Let $v_1 = F_i \cap F_j$ and $v_2 = F_i \cap F_k$ with $v,w \in
  \preds(\mqa)$ and $j < k$. If there is a $w = F_i \cap F_\ell$ with
  $i < \ell < k$ and $w \notin \preds(\mqa)$ then neither $v_1$ nor
  $v_2$ are reachable from $w$. Thus the sink $\s_i \in F_i$ of the
  2-face $F_i$ must be located on $F_i$ between $v_1$ and $v_2$, i.e.
  there is $m$ with $j < m < k$ and $\s_i = F_i \cap F_m$. But since
  $\s_i$ is reachable from $w' \defi F_i \cap F_0$, at least one of the
  vertices $v_1$ and $v_2$ is reachable from $w' \in F_0$. This is a
  contradiction to $F_0 \cap \preds(\mqa) = \emptyset$.
\end{proof}

Now we begin the proof of condition \eqref{eq:qs} for the function
$\lambda_\Qset$ (i.e. decreasing direct successor and increasing direct
predecessor for every vertex $v \in \Qset$).  
The proof is split into various lemmas. Each applying to a slightly
different situation.

\begin{lemma} \label{lem:main1-a}
  Let $v \in \Qset_1 \cup \Qset_4$ and $i$, $j$ such that $\{v\} = F_i \cap
  F_j$. and $\alpha \in \{1,2,3,4\}$ such that $\q_i \in
  \gamma_\alpha$. If $\gamma_\alpha\vert_{[q,\q_i]} \cap F_j \neq
  \emptyset$ and either $v$ and $\q_i$ are both horizontal or both
  vertical w.r.t. $F_i$ then $\lambda_\Qset$
  satisfies \eqref{eq:qs} at $v$.
\end{lemma}
\begin{proof}
  Set $\hat{\imath} \in [0,n-1]$
  \begin{equation*}
    \hat{\imath} \defi \begin{cases}
      i+1 &\text{if } \alpha \in \{1,4\} \\
      i-1 &\text{if } \alpha \in \{2,3\}. 
    \end{cases}
  \end{equation*}
  Then $\q_{\hat{\imath}}$ is the next source on $\gamma_\alpha$.
  Define $w_1 \in \gamma_\alpha\vert_{[q,\q_i]} \cap F_j$ to be the
  last such vertex on $\gamma_\alpha\vert_{[q,\q_i]}$. Choose $\Phi$
  from Lemma \ref{lem:embeddings} as follows.
  If $\alpha=1$ choose
  \begin{equation*}
    x(\,\Phi(w)\,) \defi n - \min(w) + 1
  \end{equation*}
  and
  \begin{equation*}
    y(\,\Phi(w)\,) \defi \begin{cases}
      \max(w) & \text{if } \max(\q_i) > \max(v) \\
      n - \max(w) + 1 & \text{if } \max(\q_i) < \max(v)
    \end{cases} 
  \end{equation*}
  If $\alpha=2$ choose
  \begin{equation*}
    x(\,\Phi(w)\,) \defi \max(w)
  \end{equation*}
  and
  \begin{equation*}
    y(\,\Phi(w)\,) \defi
    \begin{cases}
      \min(w)&\text{if } \min(\q_i) > \min(v)\\
      n - \min(w) + 1&\text{if } \min(\q_i) < \min(v)
    \end{cases}
  \end{equation*}
  If $\alpha=3$ choose
  \begin{equation*}
    x(\,\Phi(w)\,) \defi \min(w)
  \end{equation*}
  and
  \begin{equation*}
    y(\,\Phi(w)\,) \defi \begin{cases}
      \max(w) & \text{if } \max(\q_i) > \max(v) \\
      n - \max(w) + 1 & \text{if } \max(\q_i) < \max(v)
    \end{cases}
  \end{equation*}
  If $\alpha=4$ choose
  \begin{equation*}
    x(\,\Phi(w)\,) \defi n-\max(w)+1
  \end{equation*}
  and
  \begin{equation*}
    y(\,\Phi(w)\,)
    \defi \begin{cases}
      \min(w)&\text{if } \min(\q_i) > \min(v)\\
      n - \min(w) + 1&\text{if } \min(\q_i) < \min(v)
    \end{cases}
  \end{equation*}
%
                                %
  Then $\Phi(\gamma_\alpha\vert_{[q,\q_i]})$ is horizontal monotone by
  \eqref{eq:mainlemma:c3:phix}.  Set $w_2 \defi \q_i$ and $v_0 \defi
  v$. Then using the definition of $y(\,\Phi(w)\,)$ above
  \begin{align}
    x(\,\Phi(w_1)\,) &\in [x(\,\Phi(v_0)\,)+1, n-1] & y(\,\Phi(w_1)\,)
    &= y(\,\Phi(v_0)\,)     \\
    x(\,\Phi(w_2)\,) &= x(\,\Phi(v_0)\,) & y(\,\Phi(w_2)\,) &\in
    [y(\,\Phi(v_0)\,)+1,n-1].
  \end{align}
  Setting $B \defi [i,n-1] \times [j,n-1]$,
  $\Phi(\gamma_\alpha\vert_{[w_1,w_2]}) \subset B$ holds by the choice
  of $x(\,\Phi(w)\,)$ and using the fact that $w_1$ was chosen to be
  the `last' vertex in $\gamma_\alpha\vert_{[q,\q_i]} \cap F_j$.
  We would like to apply Corollary~\ref{cor:paths-b} with $\omega_2
  \defi \Phi(\gamma_\alpha\vert{[w_1,w_2]})$ and $\omega_1$ being the
  empty path.  Thus we have to check for any skew edges in the image of
  $\Phi$ which might cause trouble.

  A skew edge in the image of $\Phi$ causing trouble can only be an edge
  whose pre-image is incident to a vertex $w \in F_i \vert_{[\q_i,v]}$
  and whose endpoints are intermediate vertices, i.e. are incident to
  three 2-faces of $\F$. 
  Since both vertices $v$ and $\q_i$ are horizontal respectively
  vertical on $F_i$, the skew edges can only be incident to $v$ or
  $\q_i$ on $F_i\vert_{[\q_i,v]}$. Skew edges incident to $\q_i$ do not cause
  any trouble.  
  Now suppose that $v$ and $\q_i$ are both horizontal and $v$ is
  incident to three 2-faces in $\F$. Then $j = i-2$ holds and if $\q_i
  \in \gamma_3$ it would follow, that $v \in W_2$. This is a
  contradiction to $v \in \Qset_1 \cup \Qset_4$ and thus $\q_i \in \gamma_1$
  and $v \in W_4$. But then $\gamma_1\vert_{[q,\q_i]} \cap F^V_j \not=
  \emptyset$ is a contradiction. 
  If $v$ and $\q_i$ are both vertical, the same holds with $j=i+2$ and
  $\gamma_2$ and $\gamma_4$.
  Thus there are no skew edges incident to $v$.

  Using Corollary~\ref{cor:paths-b} we can conclude that every
  directed path $\gamma':\q_{\hat{\imath}} \dto v$ intersects
  $\gamma_\alpha\vert_{[w_1,\q_i]}$, since $\Phi(\q_{\hat{\imath}})
  \notin B$. But this would result in a directed cycle since
  $\q_{\hat{\imath}}$ is the next source on $\gamma_\alpha$ after
  $\q_i$. Thus such a path $\gamma'$ cannot exist and
  $\q_{\hat{\imath}} \notin \preds(v)$.
  Setting $v' \defi F_j \cap F_{\hat{\imath}}$ then $(v,v') \in A$ and
  $\q_{\hat{\imath}} \in \preds(v')$ but $\q_{\hat{\imath}} \notin
  \preds(v)$. This proves \eqref{eq:qs:q} for $v$ and $\lambda_\Qset$.
  
  Condition \eqref{eq:qs:s} follows with $\tilde{\imath} \defi i + 1$
  if $\alpha \in \{1,4\}$ and $\tilde{\imath} \defi i-1$ if $\alpha
  \in \{2,3\}$ for $\tilde{v} \defi F_j \cap F_{\tilde{\imath}}$ using
  the same Argument on $\tilde{v}$ and the edge $(\tilde{v},v) \in A$.
\end{proof}

\begin{lemma} \label{lem:main1-b}
  Let $v \in \Qset_1 \cup \Qset_4$ and set $i,j \in
  [0,n-1]$ such that $\{v\} = F_i \cap F_j$ and $\alpha \in \{1,3\}$,
  $\beta \in \{2,4\}$ such that $\q_i \in \gamma_\alpha$ and $\q_j \in
  \gamma_\beta$.  If $\gamma_\alpha\vert_{[q,\q_i]} \cap F_j =
  \emptyset$ and $\gamma_\beta\vert_{[q,\q_j]} \cap F_i = \emptyset$
  then $v$ and $\lambda_\Qset$ satisfy conditions \eqref{eq:qs}.
\end{lemma}
\begin{proof}
  Set $\hat{\imath}, \hat{\jmath} \in [0,n-1]$
  \begin{align*}
    \hat{\imath} \defi \begin{cases}
      i+1 & \text{if } \alpha = 1 \\
      i-1 & \text{if } \alpha = 3
    \end{cases} &&
    \hat{\jmath} \defi \begin{cases}
      i+1 & \text{if } \beta = 4 \\
      i-1 & \text{if } \beta = 2.
    \end{cases}
  \end{align*}
  Then $\q_{\hat{\imath}}$ is the next source on $\gamma_\alpha$ and
  $\q_{\hat{\jmath}}$ is the next source on $\gamma_\beta$.
  
  Suppose that $\gamma_\alpha\vert_{[q,\q_i]} =
  \gamma_\beta\vert_{[q,\q_j]} = \emptyset$.  Choose $\Phi$ from Lemma
  \ref{lem:embeddings} with
  \begin{align}
    x(\,\Phi(w)\,) &\defi \begin{cases}
      n - \min(w) + 1 & \text{if } \alpha = 1\\
      \min(w) & \text{if } \alpha = 3
    \end{cases} \label{eq:mainlemma:c3:phix} \\
    y(\,\Phi(w)\,) &\defi \begin{cases}
      \max(w) & \text{if } \beta = 2\\
      n - \max(w) + 1 &\text{if } \beta = 4
    \end{cases} \label{eq:mainlemma:c3:phiy}
  \end{align}
  Then $\Phi(\gamma_\alpha\vert_{[q,\q_i]})$ is horizontal monotone by
  \eqref{eq:mainlemma:c3:phix} and $\Phi(\gamma_\beta\vert_{[q,\q_j]})$
  is vertical monotone by \eqref{eq:mainlemma:c3:phiy}. Set $w_2 \defi
  \q_i$, $w_1 \defi \q_j$ and $k,\ell \defi n$, $B \defi [i,n-1] \times
  [j,n-1]$. Then $\Phi(q) \in B$ and again as in Case 1
  \begin{align}
    x(\,\Phi(w_1)\,) &\in [x(\,\Phi(v_0)\,)+1, n-1] & y(\,\Phi(w_1)\,)
    &= y(\,\Phi(v_0)\,)   \\
    x(\,\Phi(w_2)\,) &= x(\,\Phi(v_0)\,) & y(\,\Phi(w_2)\,) &\in
    [y(\,\Phi(v_0)\,)+1,n-1].
  \end{align}
                                
  Analogously to the proof of \ref{lem:main1-a} we have to show, that
  there are no skew edges in the image of $\Phi$ which might cause any
  trouble, in order to apply Corollary~\ref{cor:paths-b}.
  Skew edges in the image of $\Phi$ causing trouble can only be edges
  whose pre-image is incident to a vertex $w \in F_i \vert_{[\q_i,v]}
  \cup F_j\vert_{[\q_j,v]}$ and whose endpoints are intermediate
  vertices, i.e. are incident to three 2-faces of $\F$.
  Since $\q_i \in \gamma_1 \cup \gamma_3$ and $\q_j \in \gamma_2 \cup
  \gamma_4$, $v \in (V_{14} \cup V_{12} \cup V_{34} \cup V_{23}) \cap
  (\Qset_1 \cup \Qset_4)$ holds. But $\mqb$ is the only intermediate vertex in
  $V_{14} \cup V_{12} \cup V_{34} \cup V_{23}$. 
  Thus there are no skew edges in the image of $\Phi$ causing any
  trouble and we can apply Corollary~\ref{cor:paths-b} with $\omega_1
  \defi \Phi(\gamma_\beta\vert_{[q,\q_j]})$ and $\omega_2 \defi
  \Phi(\gamma_\alpha\vert_{[q,\q_i]})$ to conclude that the directed
  paths $\gamma':\q_{\hat{\imath}} \dto v$ and
  $\gamma'':\q_{\hat{\jmath}} \dto v$ must intersect
  $\gamma_\beta\vert_{[q,\q_j]}$ respectively
  $\gamma_\alpha\vert_{[q,\q_i]}$. Thus both cannot exist. And we can
  find a neighbor $v'$ of $v$ with $(v,v') \in A$ and $\lambda_\Qset(v') <
  \lambda_\Qset(v)$ which proves \eqref{eq:qs:q}. Again condition
  \eqref{eq:qs:s} is verified analogously.
\end{proof}

\begin{corollary} \label{cor:main1}
  $\lambda_\Qset$ satisfies conditions \eqref{eq:qs} on $\Qset_1$.
\end{corollary}
\begin{proof}
  Let $v \in \Qset_1$ and $i \defi \min(v)$, $j \defi \max(v)$. By the
  definition of $\Qset_1$ and Lemma~\ref{lem:Q-partition} $\q_i, \q_j \in
  \Qset_1$ and both $v$ and $\q_i$ are horizontal on $F_i$ and both $v$ and
  $\q_j$ are vertical on $F_j$. Thus we can apply Lemma~\ref{lem:main1-b} 
  if $\gamma_\alpha\vert_{[q,\q_i]} \cap F_j =
  \emptyset$ and $\gamma_\beta\vert_{[q,\q_j]} \cap F_i =
  \emptyset$. Otherwise we can apply Lemma~\ref{lem:main1-a}.
\end{proof}

\begin{lemma} \label{lem:main3-a}
  Let $v \in \Qset_4 \cap (V_{11} \cup V_{13})$ and $i \defi \min(v)$, $j
  \defi \max(v)$, i.e.  $v = F^H_i \cap F^V_j$ and $\q_i \in \gamma_1$.
  If $F^V_j \cap \gamma_1 = \emptyset$ then $\lambda_\Qset$ satisfies
  \eqref{eq:qs} at $v$.
\end{lemma}
\begin{proof}
  The vertex 
  $\q_{i-1}$ is the next source after $\q_i$ on $\gamma_1$.
  We will show, that every path from $\q_{i-1}$ to $v$ must cross
  $\gamma_1\vert_{[q,\q_i]} \cup \gamma_4$. Which is a contradiction,
  since $\q_{i-1} \in \succs(\q_i)$ and $\q_i \in \succs(\mqb)$.
                               
  We choose the embedding $\Phi$ among those in
  Lemma~\ref{lem:embeddings} to be
  \begin{align}
    x(\,\Phi(w)\,) &\defi n-\min(w)+1 & y(\,\Phi(w)\,) &\defi 
    n - \max(w) + 1
  \end{align}
  By Lemma~\ref{lem:embeddings} $\Phi(\gamma_1)$ is a horizontal
  monotone source-path.  By Lemma~\ref{lem:embeddings-sinkfree} there
  is a sink-free fence $D'$ such that $\Phi(\preds(\mqa))$ is a
  subgraph of $D'$.

  If $\max(\q_i) < \max(v)$, then $\gamma_1\vert_{[q,\q_i]} \cap F^V_j
  \not= \emptyset$.  Thus $\max(\q_i) > \max(v)$.
  And we construct $D'$ in the following way. First $D'$ does not have
  any more skew edges than the image of $\Phi$. Secondly
  $\Phi(\gamma_4)$ can be extended in $D'$ by a straight vertical
  source-path 
  \begin{equation*}
    \omega' = \Phi(\mqb) (\,x(\Phi(\mqb)),
    y(\Phi(\mqb))+1\,) (\,x(\Phi(\mqb)), y(\Phi(\mqb))+2\,) \ldots
    (\,x(\Phi(\mqb)), y(\Phi(v))\,),
  \end{equation*}
  each vertex (except $\Phi(\mqb)$)
  being the source of the horizontal line.  This is possible since all
  2-faces being the preimage of those vertical lines have sources in
  $\gamma_1$ and $\gamma_3$ only, and thus horizontal sources. 
  
  We would like to apply Corollary~\ref{cor:paths-b} with $\omega_1
  \defi \Phi(\gamma_4) \circ \omega'$ and $\omega_2 \defi
  \Phi(\gamma_1)$. Thus we have to check for any skew edges which may
  cause trouble.  Since $v$ and $\q_i$ are both horizontal vertices on
  $F_i$ there are no skew edges adjacent to the image of $F_i^H
  \vert_{[\q_i,v]}$.  Thus we need to check the horizontal line of
  $\Phi(v)$ between $\Phi(v)$ and the endpoint of $\omega'$.  But the
  only possible skew edges are those adjacent to vertices to $F^V_j$.
  But since $\q_j$ is a horizontal source, the skew edges adjacent to
  $\Phi(F^V_j)$ do not cause any trouble.  Thus we can apply
  Corollary~\ref{cor:paths-b}. Since $\omega'$ does not have any
  preimage, we can conclude that any path from $\q_{j-1}$ to $v$ must
  intersect either $\gamma_1$ or $\gamma_4$.
\end{proof}

\begin{corollary} \label{cor:main3-b}
  Let $v \in \Qset_4 \cap (V_{22} \cup V_{24})$ and $i \defi \min(v)$, $j
  \defi \max(v)$, i.e.  $v = F^H_i \cap F^V_j$ and $\q_j \in \gamma_2$.
  If $F^V_j \cap \gamma_2 = \emptyset$ then $\lambda_\Qset$ satisfies
  \eqref{eq:qs} at $v$.
\end{corollary}
\begin{proof}
  This corollary follows from Lemma~\ref{lem:main3-a} by reversing the
  order of the 2-faces $\F$ while keeping $F_0$ fixed. 
\end{proof}

\begin{corollary} \label{cor:main2}
  $\lambda_\Qset$ satisfies conditions \eqref{eq:qs} on $\Qset_4$.
\end{corollary}
\begin{proof}
  Let $v \in \Qset_4$ and $i \defi \min(v)$, $j \defi \max(v)$. 
  
  If $\q_i$ is horizontal and $\q_j$ is vertical we can either apply
  Lemma~\ref{lem:main1-b} or Lemma~\ref{lem:main1-a} as in
  Corollary~\ref{cor:main1}.

  If both $\q_i$ and $\q_j$ are horizontal, then by
  Lemma~\ref{lem:monotone-paths} and by the definition of $\Qset_4$ we an
  deduce that $v \in \Qset_4 \cap (V_{11} \cup V_{13})$ holds. And thus we
  can either apply Lemma~\ref{lem:main3-a} or again
  Lemma~\ref{lem:main1-a}.
  
  If both $\q_i$ and $\q_j$ are vertical, then---analogously to the
  previous case--by Lem\-ma~\ref{lem:monotone-paths} and by the
  definition of $\Qset_4$ we an deduce that $v \in \Qset_4 \cap (V_{22} \cup
  V_{24})$ holds. And thus we can either apply Lemma~\ref{cor:main3-b}
  or Lemma~\ref{lem:main1-a}.
\end{proof}

Now we prove the main proposition of this subsection, which puts
everything together.

\begin{proposition} \label{prop:main-qs}
  For any AUSO $D=(V,A)$ of a dual cyclic 4-polytope $C^\polar(n)$,
  there is a partition $\Pi_\Qset$ of $\Qset = \preds(\mqa) \cup \preds(\mqb)$
  with $\# \Pi_\Qset \le 4$ and for all $W \in \Pi_\Qset$ there is a monotone
  decreasing function $\lambda_W : W \rightarrow \Z$ satisfying
  \eqref{eq:qs} for all $v \in W$ and with $\#\lambda_W(W) \le n$ and $D
  / \Pi_\Qset$ is acyclic.
\end{proposition}

\begin{proof}
  W.l.o.g. we assume $\mqa \notin \preds(\mqb)$. Thus there are two
  cases.
  
  \emph{Case 1: $\mqb \notin \preds(\mqa)$.}  Then $\Qset_2 = \Qset_3 = \Qset_4 =
  \emptyset$. And $\lambda_\Qset$ satisfies \eqref{eq:qs} on $\Qset_1 =
  \preds(\mqa)$ by Corollary~\ref{cor:main1}. Since neither $\mqa$ nor
  $\mqb$ are sinks of 2-faces in $\F$, we can exchange $\mqa$ and
  $\mqb$ yielding that $\lambda_\Qset$ satisfies \eqref{eq:qs} on
  $\preds(\mqb)$. Thus we may set
  \begin{gather}
    \Pi_\Qset \defi \left\{ \preds(\mqa) \cap \preds(\mqb), \preds(\mqa)
      \setminus \preds(\mqb), \preds(\mqb) \setminus \preds(\mqb)
    \right\} \\
    \lambda_{\preds(\mqa) \cap \preds(\mqb)} \defi \lambda_\Qset \quad
    \lambda_{\preds(\mqa) \setminus \preds(\mqb)} \defi \lambda_\Qset \quad
    \lambda_{\preds(\mqb) \setminus \preds(\mqa)} \defi \lambda_\Qset
  \end{gather}
  Clearly $D[\Qset]/\Pi_\Qset$ is acyclic and each $\lambda_W$ is a monotone
  decreasing function satisfying \eqref{eq:qs} and $\#\lambda_W(W) \le
  n$.
    
  \emph{Case 2: $\mqb \in \preds(\mqa)$.}  Thus $\Qset=\preds(\mqa)$. And
  we use the partition together with the functions defined earlier in
  this section.
  \begin{gather}
    \Pi_\Qset \defi \left\{ \Qset_1, \Qset_2, \Qset_3, \Qset_4 \right\} \\
    \lambda_1 \defi \lambda_\Qset \quad
    \lambda_2 \defi \min \quad
    \lambda_3 \defi -\max \quad
    \lambda_4 \defi \lambda_\Qset \quad
  \end{gather}
  By Lemma~\ref{lem:Q-partition} $\Pi_\Qset$ is a partition of $\Qset$ and
  $D[\Qset]/\Pi_\Qset$ is acyclic. 
  By Lemma~\ref{lem:main-Q2} $\lambda_2$ and by
  Corollary~\ref{cor:main-Q3} the functions $\lambda_2$ and
  $\lambda_3$ are monotone decreasing and satisfy condition
  \eqref{eq:qs} on $\Qset_2$ respectively $\Qset_3$.
  By Corollary~\ref{cor:main1} and Corollary~\ref{cor:main2}
  $\lambda_1$ and $\lambda_4$ satisfy \eqref{eq:qs} on $\Qset_1$
  respectively $\Qset_4$ and clearly both functions are monotone
  decreasing.
\end{proof}

No we reformulate Proposition~\ref{prop:main-qs} into two corollaries one
dealing with the set $\Qset$ and the other one dealing with the set $\Sset$.

\begin{corollary} \label{cor:main-q}
  For any AUSO $D=(V,A)$ of a dual cyclic 4-polytope $C^\polar(n)$,
  there is a partition $\Pi_\Qset$ of $\Qset = \preds(\mqa) \cup \preds(\mqb)$
  with $\# \Pi_\Qset \le 6$ and for all $W \in \Pi_\Qset$ there is a monotone
  decreasing function $\lambda_W : W \rightarrow \Z$ with decreasing
  direct successors with respect to $W$
  and with $\#\lambda_W(W) \le n$ and $D / \Pi_\Qset$ is acyclic.
\end{corollary}

\begin{corollary} \label{cor:main-s}
  For any AUSO $D=(V,A)$ of a dual cyclic 4-polytope $C^\polar(n)$,
  there is a partition $\Pi_\Sset$ of $\Sset = \succs(\msa) \cup \succs(\msb)$
  with $\# \Pi_\Sset \le 6$ and for all $W \in \Pi_\Sset$ there is a monotone
  decreasing function $\lambda_W : W \rightarrow \Z$ with decreasing
  direct successors with respect to $W$
  and with $\#\lambda_W(W) \le n$ and $D / \Pi_\Sset$ is acyclic.
\end{corollary}

\begin{proof}
  We exploit the fact, that the inverse of an AUSO is again an AUSO.
  Going to the inverse exchanges the sets $\Qset$ and $\Sset$, the functions
  $\lambda_W$ (defined on $\Qset$) become monotone increasing (on $\Sset$),
  thus $(-\lambda_W)$ is monotone decreasing (on $\Sset$). Since each
  $\lambda_W$ satisfies \eqref{eq:qs:s} (on $\Qset$), $(-\lambda_W)$ has
  decreasing direct successors (on $\Sset$).
\end{proof}

\subsection {\boldmath Investigating the Vertex Set $R$}
\label{sect:analysis:r}

Now we want to proof a theorem similar to Corollaries \ref{cor:main-q}
and \ref{cor:main-s} for the set $R = V \setminus (\Qset \cup \Sset)$.
The order on the 2-faces in $\F$ induces an orientation on the edges
of each 2-face $F_k \in \F$. We exploit this to assign a sign to each
directed edge $e = (v,w) \in A$.
Define the sign $\sigma : A \rightarrow \{+,-\}$ of $e$ as
\begin{equation}
  \sigma(e) \defi \begin{cases}
    +& \max(w) = \max(v) + 1 \mod n 
    \text{ or } \min(w) = \min(v) +1 \mod n\\
    -& \max(w) = \max(v) - 1 \mod n 
    \text{ or } \min(w) = \min(v) - 1 \mod n
  \end{cases}.
\end{equation}
Note that this definition does not depend on choosing a particular
2-face as $F_0$.
Since there are no sources or sinks of the 2-faces $\F$ in the vertex
set $R$, we can define the following sign function
$\sigma_F:R\rightarrow \{+,-\}$ for vertices $v \in R$ incident to a
given 2-face $F$.
\begin{equation}
  \sigma_F(v) \defi \begin{cases}
    +& \text{For all $e \in E(F)$ with $v \in e$, $\sigma(e) = +$ holds.}\\
    -& \text{For all $e \in E(F)$ with $v \in e$, $\sigma(e) = -$ holds.}\\
  \end{cases}
\end{equation}
Now we define the following four subsets of $R$.
\begin{align*}
  R^{++} \defi & \left\{ v \in R \,:\, \sigma_{F_{\min(v)}}(v) = +
    \text{ and }\sigma_{F_{\max(v)}}(v) = + \right\} \\
  R^{--} \defi & \left\{ v \in R \,:\, \sigma_{F_{\min(v)}}(v) = -
    \text{ and }\sigma_{F_{\max(v)}}(v) = - \right\} \\
  R^{-+} \defi & \left\{ v \in R \,:\, \sigma_{F_{\min(v)}}(v) = -
    \text{ and }\sigma_{F_{\max(v)}}(v) = + \right\} \\
  R^{+-} \defi & \left\{ v \in R \,:\, \sigma_{F_{\min(v)}}(v) = +
    \text{ and }\sigma_{F_{\max(v)}}(v) = - \right\} \\
\end{align*}

Our investigation will continue by showing,
that these four sets are disconnected (i.e. are separated by $\Qset$ and
$\Sset$), and that thus we can easily find effectively decreasing
functions.

\begin{lemma} \label{lem:sign-transportation}
  Let $e = (v,w)$ be an arbitrary edge with $v,w \in R$.  Let
  $(i,j,k)$ be a triple of indices such that $v \in F_i \cap F_k$, $w
  \in F_j \cap F_k$, $e \in F_k$ and $F_i$ and $F_j$ are neighbors
  (i.e. $i-j \mod n =1$).
  
  Then $\sigma_{F_i}(e_i) = \sigma_{F_j}(e_j)$, for all edges $e_i \in
  F_i$ incident to $v$ and all edges $e_j \in F_j$ being incident to
  $w$.
\end{lemma}

\begin{figure}
  \begin{center}
    \begin{overpic}[tics=5,scale=.4]{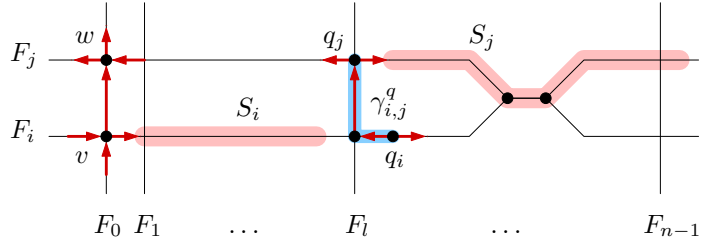}
      \put(11.5,0){$F_0$}
      \put(17,0){$F_1$}
      \put(30,0){$\ldots$}
      \put(46,0){$F_l$}
      \put(66,0){$\ldots$}
      \put(87,0){$F_{n-1}$}
      
      \put(0,23.5){$F_j$}
      \put(0,13){$F_i$}

      \put(9,9.5){$v$}
      \put(9,26){$w$}
      \put(43,26){$\q_j$}
      \put(49.5,17){$\gamma^\q_{i,j}$}
      \put(51.5,9.5){$\q_i$}

      \put(31,16){$\Sset_i$}
      \put(63,26){$\Sset_j$}
    \end{overpic}
  \end{center}
  \caption{ \label{fig:sign-transportation}
    Illustrating the proof of Lemma~\ref{lem:sign-transportation}.
}
\end{figure}

\begin{proof}
  Since $\sigma$ does not depend on choosing the 2-face $F_0$, we can
  assume w.l.o.g. $k=0$.
  Let $e'$ be the last edge of the directed path $\gamma^\q_{ij}$ and
  let $\ell$ be the index such that $e' \in F_\ell$. By the special
  structure of the paths $\gamma^\q_{ij}$ (constructed in
  Section~\ref{sect:prelims:dcpausos}), $e'$ is incident to $F_i$ and
  $F_j$.
  Let $s,t \in [0, n-1]$ such that $\{\q_i\} = F_i \cap F_s$ and
  $\{\q_j\} = F_j \cap F_t$.
  Suppose that $\sigma_{F_i}(v) = +$ and $\sigma_{F_j}(w) = -$.  
  Thus we can conclude that $\s_i$ and $\s_j$ are located in the
  following segments of $F_i$ respectively $F_j$.
  \begin{align} \label{eq:lem:sign-transportation:eq1}
    \s_i \in F_i \cap \left( F_1 \cup F_2 \cup \cdots \cup
      F_{s-1}\right) 
    && \s_j \in F_j \cap \left( F_{t+1} \cup F_{t+2}
      \cup \cdots \cup F_{n-1}\right)
  \end{align}
  $v,w \notin \gamma^\q_{i,j}$, since $\gamma^\q_{i,j} \subset \Qset$. Thus
  \eqref{eq:lem:sign-transportation:eq1} implies that $\s_i$ and $\s_j$
  are even located in the following segments of $F_i$ respectively
  $F_j$.
  \begin{align*}
    \s_i &\in \Sset_i \defi F_i \cap \left( F_1 \cup F_2 \cup \cdots \cup
      F_{\ell-1}\right) \\
    \s_j &\in \Sset_j \defi F_j \cap \left( F_{\ell+1} \cup F_{\ell+2}
      \cup \cdots \cup F_{n-1}\right)
  \end{align*}
  But thus $\gamma^\s_{ij}$ must intersect either $F_0$ or $F_\ell$.
  The first case violates $v,w \in R$ the latter violates
  Corollary~\ref{cor:qs-paths}.
\end{proof}

\begin{lemma} \label{lem:separation-R}
  For all pairs of vertices $x,y \in R$ being contained in different
  sets among the four sets $R^{++}$, $R^{--}$, $R^{-+}$, and $R^{-+}$,
  there is no path connecting $x$ and $y$ only using vertices in $R$.
\end{lemma}
\begin{proof}
  For every edge $e = (v,w)$ there is a triple $(i,j,k)$ of indices
  such that $v \in F_i \cap F_k$, $w \in F_j \cap F_k$, $e \in F_k$
  and $F_i$ and $F_j$ are neighbors (i.e. $i-j \mod n =1$), Lemma
  \ref{lem:sign-transportation} implies, that for all pairs $x,y \in
  R$ joined by a path though $R$ the endpoints $x$ and $y$ must be
  contained in the same set $R^{++}$, $R^{--}$, $R^{-+}$, or $R^{-+}$.
\end{proof}

\begin{proposition} \label{prop:main-R} For any AUSO $D=(V,A)$ of a dual
  cyclic 4-polytope $C^\polar(n)$, there is a partition $\Pi_R$ of $R$
  with $\# \Pi_R \le 4$ and for all $W \in \Pi_R$ there is a
  effectively decreasing function (with respect to $W$) $\lambda_W : W
  \rightarrow \Z$
  and with $\#\lambda_W(W) \le n$, and $D[R] / (\Pi_R \cup \{\Qset,\Sset\})$ is
  acyclic.
\end{proposition}
\begin{proof}
  Clearly we set $\Pi_R = \{R^{++},R^{--},R^{-+},R^{+-}\}$ and we
  define the following functions
  \begin{align*}
    \lambda_{++}& : R^{++} \rightarrow \Z, \lambda_{++}(v) \defi -\min(v) \\
    \lambda_{--}& : R^{--} \rightarrow \Z, \lambda_{--}(v) \defi \min(v) \\
    \lambda_{-+}& : R^{-+} \rightarrow \Z, \lambda_{-+}(v) \defi \min(v) \\
    \lambda_{+-}& : R^{+-} \rightarrow \Z, \lambda_{+-}(v) \defi -\min(v) \\
  \end{align*}
  The functions are effectively decreasing (with respect to each
  $R^{**}$).  The acyclicity of $D/ \Pi_R$ follows from Lemma
  \ref{lem:separation-R} and the fact that the set $\Qset$ and $\Sset$
  is a source- respectively sink-set.
\end{proof}

\subsection{Proving the Main Theorem}
The preceeding analysis yields the following proof of the Main Theorem~\ref{thm:main}.

\begin{proof}
  We use Corollaries~\ref{cor:main-q} and \ref{cor:main-s}, and
  Proposition~\ref{prop:main-R};
  the partition 
\begin{equation*}
  \Pi \defi \{ \Qset_1, \Qset_2, \Qset_3, \Qset_4, 
  \Sset_1, \Sset_2, \Sset_3, \Sset_4, R^{++}, R^{--}, R^{-+}, R^{+-}\}
\end{equation*}
   and
  the functions defined in the above theorems.
  Then we have a partition $\Pi$ of the vertices $V$ of constant size
  $\# \Pi \le 12$. And for every $W \in \Pi$ there is an effectively
  decreasing function (with respect to monotone) $\lambda_W : W
  \rightarrow \Z$ and with $\#\lambda_W(W) \le n$, and $D / \Pi$ is
  acyclic.

  Thus we can apply Theorem~\ref{thm:runtime} which proves, that
  random edge does not use more then $O(12 \cdot n) = O(n)$ pivot
  steps starting at an arbitrary vertex.
\end{proof}

\section {Remarks}
\label{sect:rem}

Now that we have seen that random edge takes only linear expected
running time on dual cyclic 4-polytopes, we would like to shed some
light on the question whether this is an inherent property of random
edge or it is rather caused by the very specific structure of the
considered dual cyclic 4-polytopes.

\subsection {Products of Two Polygons}
\label{sect:polygons}

Let $C_n$ denote the (regular) $n$-gon in the plane.  Let $P$ be a
4-polytope which is combinatorially equivalent to $C_n \times C_m$,
where $\times$ denotes the usual product of sets.
Then $P$ is called a (combinatorial) product of two polygons.

\begin{figure}
  \begin{center}
    \includegraphics[width=.5\textwidth]{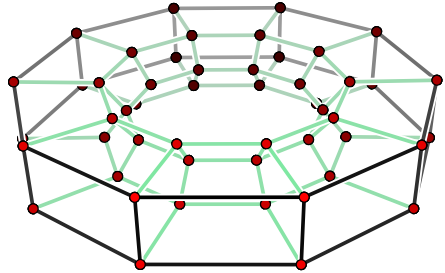}
  \end{center}
  \caption{ \label{fig:10x5} %
    The Schlegel diagram of a $10\text{-gon} \times 5\text{-gon}$.
    The (large) 2-faces are the 10-gons (horizontal) and the 5-gons
    (vertical).
  }
\end{figure}

Such polytopes have a very nice (and simple) combinatorial structure,
which resembles that of dual cyclic polytopes in some important ways.
Consider the polytope $P$.  It has two sets of large 2-faces defined
by the ``one vertex of the one polygon $\times$ the entire other
polygon''.  Each of these two sets comes with a neighborhood structure,
since the 3-faces (facets) are all defined as ``one edge of the one
polygon $\times$ the entire other polygon''.  %
(Compare to Figure~\ref{fig:10x5}.)
Thus we can achieve the same results as in
Section~\ref{sect:prelims:dcpausos}.  And we can apply the whole
machinery of Section~\ref{sect:analysis} to show that random edge
takes only $O(n+m)$ expected number of steps on $P$.  
In fact the combinatorial structure of products of two polygons is
simpler than that of dual cyclic 4-polytopes and the analysis can be
simplified and becomes considerably shorter.

\subsection {Random Facet on Dual Cyclic 4-Polytopes}
\label{sect:rf}

As mentioned in Section~\ref{sect:intro} random facet was the first
randomized pivot rule for which a sub-exponential upper bound on the
expected number of steps was proven (see~\cite{Kalai92}
and~\cite{MSW96}).
We will show that there are AUSOs on dual cyclic 4-polytopes, such
that random facet will visit at least $\Omega(n^2)$ vertices starting
at the global source.
Thus proves that for dual cyclic 4-polytopes, random edge is provably
faster than random facet.

There are several variants of the random facet rule, which differ on
how to proceed at 1-vertices (sinks of the facets).  Here we will
stick to the following definition of random facet taken from
G\"artner, Henk \& Ziegler \cite[p.~350]{GHZ98}.
See~\cite{KMSZ05} for a comparison with Kalai's original rule
in~\cite[p.~228]{Kalai97} and also a variant from
G\"artner~\cite{Gartner98}.
\begin{quote}
At each non-optimal vertex $v$ follow the (unique) outgoing edge if
$v$ is a 1-vertex. Otherwise choose one facet $f$ uniformly at random
containing $v$ and 
solve the problem restricted to $f$ by applying random facet
recursively. 
\end{quote}

The construction yields the same result for the other definitions of
random facet.  We will stick to this one, since it follows paths of
1-vertices deterministically, making the analysis slightly simpler.
It uses blocks of twelve large 2-faces. %
Let $k$ be the number of blocks used then we consider the polytopes
$C^\polar(n)$ with $n = 12 k + 1$.  The extra facet is needed to make
the construction an AUSO.

\begin{figure}
  \begin{center}
    \begin{overpic}[tics=5,scale=.7]{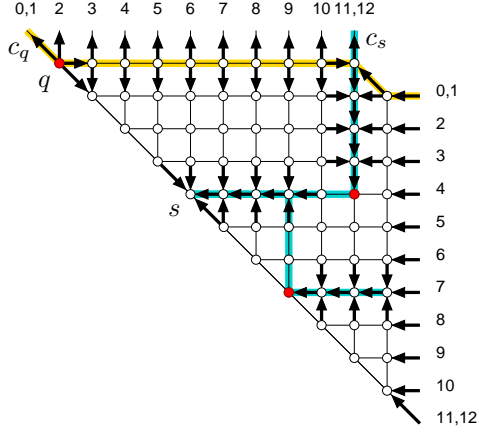}
      \put(33,45.5){$s$}
      \put(75.5,82.5){$c_s$}
      \put(-1.5,80.5){$c_q$}
      \put(5,73.5){$q$}
    \end{overpic}
  \end{center}
  \caption{ \label{fig:dcp-rf-small} %
    The constructed AUSO for $k=1$ block. 
    The (re)starting vertices of random facet are indicated.
  }
\end{figure}

Figure~\ref{fig:dcp-rf-small} depicts the constructed AUSO for $k=1$
block which we call $P_1$.  
To keep the picture simple, only the sources and sinks of all large
2-faces are indicated by oriented edges.  This fixes the orientation
of all other edges.  All vertices at which random facet may restart
are indicated.  We will call these vertices the restarting vertices.

\begin{figure}
  \begin{center}
    \begin{overpic}[tics=5,scale=.7]{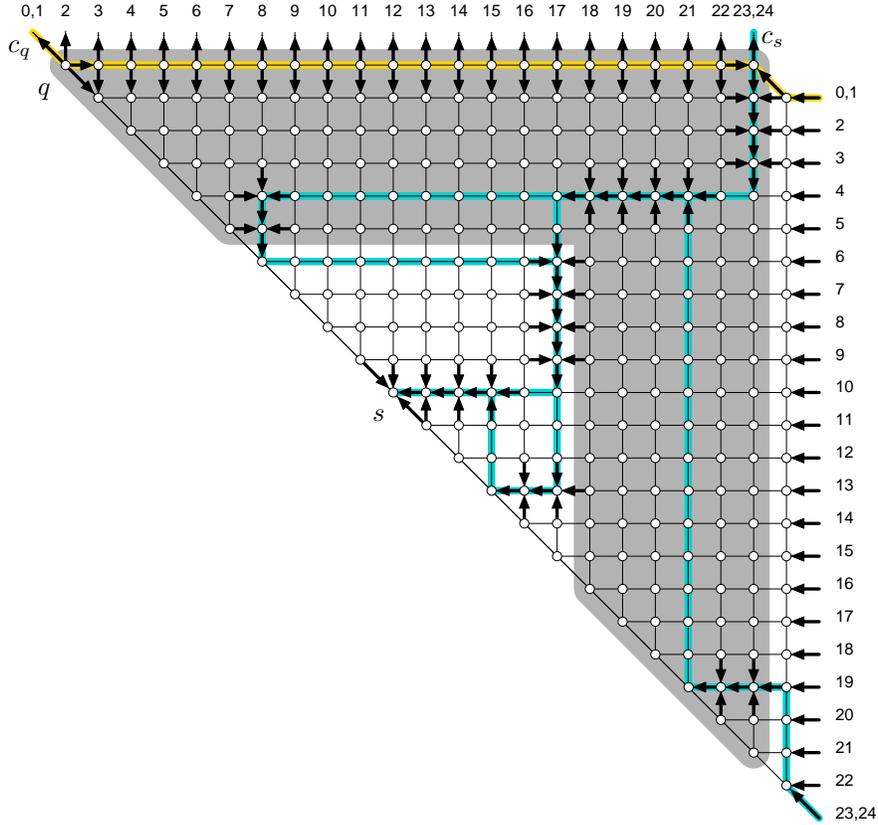}
      \put(41,47){$s$}
      \put(86.3,91){$c_s$}
      \put(-1.5,90){$c_q$}
      \put(2,85){$q$}
    \end{overpic}
  \end{center}
  \caption{ \label{fig:dcp-rf} %
    The constructed AUSO for $k=2$ blocks. 
    The (re)starting vertices of random facet are indicated.
    The shaded area indicates the 2-faces of the first block. %
  }
\end{figure}

Figure~\ref{fig:dcp-rf} depicts the constructed AUSO for $k=2$ blocks
which we call $P_2$.
The shaded area indicates the 2-faces of the first block.  
To add a new block, the new twelve 2-faces are added in the middle of
the 2-faces, i.e. in our case to get from $P_1$ to $P_2$ we added
twelve 2-faces between $F_5$ and $F_6$ in $P_1$.
Thus the global sink of $P_1$ becomes the sinks of the three 2-faces
$F_4$, $F_5$, and $F_{18}$, with coordinates $(4,18)$, $(5,18)$, and
$(4,18)$.  Now we shift the sinks of the 2-faces $F_4$ and $F_5$ by
ten coordinates/2-faces to the left, i.e. they have now coordinates
$(4,8)$ and $(5,8)$.  The resulting AUSO is $P_2$.

We will show that the shortest path possibly taken by random facet
starting at the global source $q$ contains more than $\Omega(n^2)$
vertices.
First check that random facet is restarted at the indicated vertices
only.  And that it moves from an inner restarting vertex to one of the
next restarting vertices on the diagonal.  And from an outer
restarting vertex it moves to the next interior one.
Thus we can move from Block $i$ to the next block $i+1$ only.  Further
more, when entering a new block, random facet needs at least $(k-i) 12
- 2$ steps, i.e. almost twelve steps for each coming block.
This results in an overall running time of at least
\begin{equation*}
  \sum_{i=1}^k \left( (n - i) 12  -2 \right)%
  = \sum_{j=1}^{k-1} \left( (n - i) 12  -2 \right)%
  = \frac {(k-1)(k-2)}{2} - 2(k-1) %
  = \Omega(n^2)
\end{equation*}
Thus we have proven the following theorem.
\begin{theorem} \label{thm:rf}
  There are AUSOs of $C^\polar(n)$ for $n = 12 k + 1$ such that random
  facet takes at least $\Omega(n^2)$ steps.
\end{theorem}  

Note that we have bounded the length of the shortest path possibly
taken by random facet.  Thus our lower bound holds for any random
choices and not just for the expected number of steps.
Theorem~\ref{thm:rf} even holds for any recursive pivot-rule
proceeding via incident facets.

\subsection {Conclusion}

Despite the very specific structure of dual cyclic 4-poly\-topes we were
able to separate random edge and random facet.  Similar combinatorial
properties can be found in other 4-polytopes like the product of two
polygons.  Thus maybe the results presented in this paper can be
extended to a broader class of 4-polytopes.

Nevertheless any approach to analyze random edge using only the
combinatorial notion of AUSOs must fail to give good upper bounds for
large dimensions due to the lower bounds given by Matou\v{s}ek and
Szab\'{o} in~\cite{MS04}.
Thus--as for random facet--more geometry is needed to beat the
exponential lower bounds.  One way to find more geometric properties
might be to develop further ideas for small dimensions.

\bibliographystyle{amsplain}
\bibliography{re_on_dcp}

\end{document}